\documentclass[notitlepage,leqno,10pt]{article}

\usepackage{latexsym}
\usepackage{amsmath}
\usepackage{amsthm}
\usepackage{amsfonts}
\usepackage{amssymb}
\usepackage{graphicx}

\renewcommand{\d}{\delta }
\newcommand{\D }{D }

\newcommand{\e }{\epsilon }

\newcommand{\rh }{\rho }

\newcommand{\Sig }{\Sigma}

\renewcommand{\th }{\theta }

\newcommand{\intbar}{\mathop{\int\makebox(-13.5,0){\rule[4pt]{.7em}{0.3pt}}%
\kern-6pt}\nolimits}

\newcommand{\be}{\begin{equation}}
\newcommand{\ee}{\end{equation}}
\newenvironment{pf}{\noindent{\sc Proof}.\enspace}{\rule{2mm}{2mm}\medskip}

\newcommand{\R}{\mathbb{R}}

\newcommand{\Rtre}{\mathbb{R}^3}
\newcommand{\Rdue}{\mathbb{R}^2}

\newcommand{\N}{\mathbb{N}}

\newcommand{\Sp}{\mathbb{S}}

\DeclareMathOperator{\graph}{graph}

\DeclareMathOperator{\spt}{spt}
\DeclareMathOperator{\diam}{diam}
\DeclareMathOperator{\gen}{genus}

\DeclareMathOperator{\B}{\mathcal{B}_\varepsilon}

\begin{document}

%\author{{\sc Andrea Mondino, Johannes Schygulla}}
\author{Andrea MONDINO$^{1}$, Johannes SCHYGULLA$^{2}$}

\date{}

\title{Existence of immersed spheres minimizing curvature 
functionals in non-compact $3$-manifolds}

\newtheorem{lem}{Lemma}[section]
\newtheorem{pro}[lem]{Proposition}
\newtheorem{thm}[lem]{Theorem}
\newtheorem{rem}[lem]{Remark}
\newtheorem{cor}[lem]{Corollary}
\newtheorem{df}[lem]{Definition}
\newtheorem*{Theorem}{Theorem}
\newtheorem*{Lemma}{Lemma}
\newtheorem*{Proposition}{Proposition}
\newtheorem*{claim}{Claim}
\maketitle

\footnotetext[1]{Scuola Normale Superiore, Piazza dei Cavalieri 7, 56126 Pisa, Italy, E-mail address: andrea.mondino@sns.it. \\A. Mondino is supported by a Post Doctoral Fellowship in the  ERC grant "Geometric Measure Theory in non Euclidean Spaces" directed by Prof. Ambrosio.} 
\footnotetext[2]{Mathematisches Institut, Universit\"at Freiburg, Eckerstra{\ss}e 1, 79104 Freiburg, Germany, E-mail address: johannes.schygulla@math.uni-freiburg.de}

\begin{abstract}
We study curvature functionals for immersed $2$-spheres 
in non-compact, three-dimensional Riemannian manifold $(M,h)$ without boundary. First, under the assumption that $(M,h)$ is the euclidean 3-space endowed with a semi-perturbed metric with perturbation small in $C^1$ norm and of compact support, we prove that if there is some point $\overline{x} \in M$ with scalar curvature $R^M(\overline{x})>0$ then there exists a smooth embedding  $f:\Sp^2 \hookrightarrow M$ minimizing the Willmore functional $\frac{1}{4}\int |H|^2$, where $H$ is the mean curvature. Second, assuming that $(M,h)$ is of bounded geometry (i.e. bounded sectional curvature and strictly positive injectivity radius) and asymptotically euclidean or hyperbolic we prove that if there is some point $\overline{x} \in M$ with scalar curvature $R^M(\overline{x})>6$ then there exists a smooth immersion  $f:\Sp^2 \hookrightarrow M$ minimizing the functional  $\int (\frac{1}{2}|A|^2+1)$, where $A$ is the second fundamental form. Finally, adding the bound $K^M \leq 2$ to the last assumptions,  we obtain a smooth minimizer $f:\Sp^2 \hookrightarrow M$ for the functional 
$\int (\frac{1}{4}|H|^2+1)$. The assumptions of the last two theorems are satisfied in a large class of $3$-manifolds arising as spacelike timeslices solutions of the Einstein vacuum equation in case of null or negative cosmological constant.  
\end{abstract}
\bigskip\bigskip

\begin{center}
\bigskip

\centerline{\bf AMS subject classification: }
53C21, 53C42, 58E99, 35J60

\end{center}

\section{Introduction}\label{s:in}
The present work follows the paper \cite{KMS} by Kuwert and the authors about the minimization of curvature functionals in Riemannian 3-manifolds under global conditions on the curvature of the ambient space. The aforementioned work is focalized in the case the ambient 3-manifold is compact and develop existence and regularity theory taking inspiration from \cite{SiL}. The present paper instead is concerned about the non-compact situation and relies on the regularity theory developed there. Let us point out that the study of curvature functionals, in particular of the Willmore functional, in the euclidean flat space is a topic of great interest in the contemporary research (see for instance the papers of Li \& Yau \cite{LY}, Kuwert \& Sch\"atzle \cite{KS}, Rivi\'ere \cite{Riv}, Simon \cite{SiL}, etc.); the previous \cite{KMS} and the present work are an attempt to open the almost unexplored field of the corresponding problems in non constantly curved Riemannian 3-manifolds under global geometric conditions.

Here we consider essentially two problems: first we minimize the Willmore functional among immersed spheres in $\Rtre$ endowed with a semi-perturbed metric; second we minimize related curvature functionals in non-compact Riemannian $3$-manifolds under global and asymptotic conditions on the metric. As we will remark later in the Introduction the assumptions will include a large class of manifolds naturally arising in General Relativity. Let us start discussing the first problem.

Let $h=h_{\mu\nu}$ be a symmetric bilinear form in $\Rtre$ with compact support. Denote by
$$\|h\|_{C^0}:= \sup_{x\in \Rtre} \sup_{u,v\in S^2} |h(x)(u,v)|, \quad \|Dh\|_{C^0}:= \sup_{x\in \Rtre} \sup_{u,v,w\in S^2} |D_w (h(x)(u,v))|, $$
where $D_w$ is just the directional derivative, and let $\|h\|_{C^1}=\|h\|_{C^0}+\|Dh\|_{C^0}$.\\
\\
Consider $\Rtre$ equipped with the perturbed metric $\delta+h$, where $\delta=\delta_{\mu\nu}$ is  the standard euclidean metric. For any immersed closed surface $f:\Sigma \hookrightarrow \Rtre$ with induced metric $g = f^\ast (\delta+h)$, we consider the Willmore functional
\begin{equation} \label{eq:W}
W(f) = \frac{1}{4} \int_\Sigma |H|^2\,d\mu_g,
\end{equation}
where $H$ is the mean curvature vector.\\
\\
The first problem we study is the minimization of $W(f)$ in the class of immersed spheres in the Riemannian manifold $(\Rtre,\delta+h)$ and prove the following existence result.
\begin{thm}\label{thm:ExW}
Assume $\|h\|_{C^0}\leq \eta$ and $\|Dh\|_{C^0}\leq \theta$, and that  $\spt h\subset B^e_{r_0}(x_0)$ where $B^e_{r_0}(x_0)$ is the ball in euclidean metric of center $x_0\in \Rtre$ and radius $r_0>0$ . 
On the class $[\Sp^2,(\Rtre,\delta+h)]$ of smooth immersions $f:\Sp^2 \to (\Rtre,\delta+h)$, consider the Willmore functional 
$$
W: [\Sp^2,(\Rtre,\delta+h)] \to \R,\,W(f) = \frac{1}{4} \int_\Sigma |H|^2\,d\mu_g.
$$
Assume that the scalar curvature $R_h$ of $(\Rtre,\delta+h)$ is strictly positive in some point $\overline{x} \in \Rtre$, namely $R_h(\overline{x}) > 0$. Then for $\eta$ and $r_0\theta$ sufficiently small there exists a minimizer $f$ in $[\Sp^2,(\Rtre,\delta+h)]$ for $W$, which is actually an embedding.
\end{thm}

In asymptotically flat 3-manifolds, spheres which are critical points of related curvature functionals
have been constructed recently by the first author [Mon1, Mon2]; Lamm, Metzger \& Schulze [LMS],
see also [LM], studied instead the existence of spheres which are critical points of curvature functionals under constraints. They obtain the solutions as perturbations of round spheres using implicit function type
arguments.  

Among the aforementioned papers, the most related to the present work is \cite{Mon1}; the main difference here (beside the fact that the proofs are completely different, in the former  the author used techniques of nonlinear analysis, here we use techniques of geometric measure theory) is that in the former the perturbed metric was $C^\infty$ infinitesimally close to the euclidean metric, then with infinitesimal curvature. Here instead $\delta+h$ is assumed to be close to the euclidean metric $\delta$ just in $C^0$ norm; indeed, in  order to have $r_0\theta$ small, $\|Dh\|_{C^0}$ can be large if the support of $h$ is contained in a small ball. Moreover no restrictions are imposed on the derivatives of $h$ of order higher than one, so the Riemann curvature tensor of $(\Rtre,\delta+h)$ can be arbitrarily large. For instance, if $h_{\mu\nu}(x)=h_0(x) \delta_{\mu\nu}$ for a certain  function $h_0\in C^\infty_c(\Rtre)$, then the perturbed metric $\delta_{\mu\nu}+h_{\mu\nu}=(1+h_0)\delta_{\mu\nu}$ is conformal to the euclidean metric and a direct computation shows that $R_h=2\frac{\triangle h_0}{(1+h_0)^2}-\frac{5}{2} \frac{|dh_0|^2}{(1+h_0)^3}$; therefore  taking $h_0$ with small $C^1$ norm but with large laplacian gives a metric with arbitrarily large curvature which fits in the assumptions of Theorem \ref{thm:ExW} (notice that this example is not trivial since the Willmore functional is invariant under conformal \emph{transformations} of $\Rtre$ but \emph{not} under \emph{conformal changes of metric}).
\newline
  
The second problem we study is the minimization of Willmore-type functionals in asymptotically euclidean (or asymptotically hyperbolic) Riemannian 3-manifolds. For that let $(M,h)$ be a non-compact Riemannian 3-manifold without boundary of bounded geometry, i.e.:

$i)$ $(M,h)$ has  bounded sectional curvature:
  \begin{equation}\label{Cond:BSC}
  |K^M| \leq \Lambda<\infty
  \end{equation} 

$ii)$ $(M,h)$ has strictly positive injectivity radius:
 \begin{equation}\label{Cond:InjRad}
 Inj(M,h) \geq \bar{\rho}>0.
 \end{equation}
  We assume that either

$iiia)$   $(M,h)$ is asymptotically euclidean in the following very general sense: there exist compact subsets $\Omega_1 \subset \subset M $ and $\Omega_2 \subset \subset \Rtre$ such that
\begin{equation}\label{eq:defAsFlat}
(M \backslash \Omega_1) \text{ is isometric to } (\Rtre\backslash \Omega_2, \d+o_1(1)),
\end{equation}
where $(\Rtre, \d+o_1(1))$ denotes the Riemannian manifold $\Rtre$ endowed with the euclidean metric $\delta_{\mu\nu}+o_1(1)_{\mu\nu}$ and $o_1(1)$ denotes a symmetric bilinear form which goes to 0 with its first derivatives at infinity, namely
$$\lim_{|x|\to \infty} (|o_1(1)(x)|+|\nabla o_1(1)(x)|)=0, \quad \text{or}$$

$iiib)$ $(M,h)$ is hyperbolic outside a compact subset, namely there exists $\Omega \subset \subset M$ such that the sectional curvature $K^M\leq 0$ on $M\setminus \Omega$.

For any immersed closed  surface $f:\Sigma \hookrightarrow M$ with induced metric 
$g = f^\ast h$ and second fundamental form $A$, we consider the  functional 
\begin{equation}\label{eq:defE1}
E_1(f):= \int_{\Sigma} \left( \frac{|A|^2}{2} +1\right ) d\mu_g,
\end{equation}
and we prove the following existence result.
\begin{thm}\label{thm:ExE1}
Let $(M,h)$ be a non compact Riemannian 3-manifold satisfying $i)$, $ii)$ and either $iiia)$ or $iiib)$ above. On the class $[\Sp^2,M]$ of smooth immersions $f:\Sp^2 \hookrightarrow M$, consider the functional 
$$
E_1: [\Sp^2,M] \to \R,\,E_1(f) = \int_{\Sp^2} \left( \frac{|A|^2}{2} +1\right ) d\mu_g.
$$
If the scalar curvature $R^M(\overline{x}) > 6$ for some point $\overline{x} \in M$, then there exists a smooth minimizer $f$ in $[\Sp^2,M]$ for $E_1$.
\end{thm}
Finally we will also discuss the following variant of Theorem \ref{thm:ExE1}.
\begin{thm}\label{thm:ExW1}
Let $(M,h)$ be a non compact Riemannian 3-manifold satisfying $i)$, $ii)$ and either $iiia)$ or $iiib)$ above. On the class $[\Sp^2,M]$ of smooth immersions $f:\Sp^2 \hookrightarrow M$, consider the functional 
$$
W_1: [\Sp^2,M] \to \R,\,W_1(f) = \int_{\Sp^2} \left( \frac{|H|^2}{4} +1\right ) d\mu_g.
$$
If the sectional curvature $K^M\le2$ and moreover the scalar curvature $R^M(\overline{x}) > 6$ for some point $\overline{x} \in M$, then there exists a smooth minimizer $f$ in $[\Sp^2,M]$ for $W_1$.
\end{thm}

\begin{rem}\label{rem:WeEucl}
Observe that if the ambient manifold $(M,h)$ is the euclidean space $(\Rtre,\d)$, then for every smooth immersion of a sphere $f:\Sp^2 \hookrightarrow \Rtre$  one has  $E_1(f)\geq W_1(f) > W(f)\geq 4\pi$. Moreover taking the sequence of round spheres $S_p^{1/n}$ of center $p$ and radius $1/n$ one gets $E_1(S_p^{1/n})=W_1(S_p^{1/n})=4\pi+\frac{4\pi }{n^2} \downarrow 4\pi$. So in the euclidean space the infimum of $W_1$ and $E_1$ is $4\pi$ and is never attained. Therefore the curvature assumptions are essentials for having the existence of a minimizer. 
\end{rem}

Before passing to an overview of the paper let us comment on the assumptions of Theorems \ref{thm:ExE1} and \ref{thm:ExW1}; we point out that a large class of 3-manifolds arising in General Relativity as spacelike timeslices
of solutions to the Einstein vacuum equation perfectly fit in our framework.  \\First of all observe that the asymptotic assumption $iiia)$ is very mild, indeed we are asking just an asymptotic $C^1$ closeness of the metric $h$ of the manifold with the euclidean metric; as explained above, this allows a lot of freedom to the curvature of $h$ which, for instance, is not constricted to vanish at infinity. Notice moreover that asymptotically spatial Schwarzschild 3-manifolds with mass (for the definition see, for instance, \cite{LMS} page 3), or the metric considered by Schoen and Yau in \cite{SY} in the proof of the Positive Mass Theorem, outside a ball centered in the origin, easily satisfy $iiia)$. 
\\Also assumption $iiib)$ is natural in General Relativity, indeed metrics which are asymptotic
to Anti–-de Sitter-–Schwarzschild metrics with  mass easily fit in $iiib)$ (for the definition see for instance \cite{NT} page 2, for the computation of the curvature see Lemma 3.1 of the same paper).
\\Therefore assumptions $iiia)$ and $iiib)$ correspond respectively to null and negative cosmological constant in the Einstein vacuum equations.
\newline

We conclude the Introduction by briefly outlining the contents of the present work. The technique adopted in the paper is the direct method in the calculus of variations, as in \cite{SiL} and \cite{KMS}: we consider a minimizing sequence of smooth immersions $\{f_k\}_{k \in \N}\subset  [\Sp^2,M]$ for the desired functional, we prove that the sequence is compact in a weak sense and does not degenerate, so there exists a weak minimizer and finally one gets the existence of a smooth minimizer by proving regularity. The main difficulty here is that in all the considered problems the ambient manifold is non compact, so a priori the minimizing sequence can become larger and larger in area and diameter, or may escape to infinity. Moreover, as in \cite{KMS}, the minimizing sequence can degenerate collapsing to a point. In order to prevent the aforementioned bad behaviors, we prove local and global estimates using the assumptions on the curvature of the ambient manifold. Then the weak compactness and the regularity follow as in \cite{KMS}.

More precisely in Section 2 we prove Theorem \ref{thm:ExW}; for that we first derive estimates on the geometric quantities in perturbed metric, then with a blow down procedure we get that the minimizing sequences stay in a compact subset and have bounded area, finally we prevent degeneration and  we apply similar methods and techniques developed by Simon in \cite{SiL} or Kuwert/Mondino/Schygulla in \cite{KMS} to conclude with Theorem \ref{thm:ExW}. 
\\In Section 3 we prove both Theorems \ref{thm:ExE1} and \ref{thm:ExW1}; for that we first show that minimizing sequences for the considered functionals, although the ambient manifold is non compact, stay in a compact subset of $(M,h)$ and do not degenerate. This enables us to apply the existence proof of \cite{KMS} and to conclude existence of minimizers for the functionals $E_1$ and $W_1$.

\begin{center}

{\bf Acknowledgments}

The first author would like to thank Prof.\,P.\,Caldiroli for proposing to study the 
Willmore functional in a semi-perturbative setting; he also would like to express his gratitude towards the supervisor of his Ph.D. Thesis Prof.\,A. \,Malchiodi and the supervisor of his Post Doctoral Fellowship Prof.\,L.\,Ambrosio for their constant support. All authors thanks Prof.\,E.\,Kuwert for fruitful conversations about the topics of the paper and acknowledge the support by the DFG Collaborative Research Center SFB/Transregio 71 and of M.U.R.S.T, within the project B-IDEAS ''Analysis and Beyond'', making possible our cooperation with mutual visits at Freiburg and Trieste. 
\end{center}

\noindent

\section{Proof of Theorem \ref{thm:ExW}}
\subsection{Geometric estimates and a monotonicity formula in perturbed setting}\label{Sec:MF}
The goal of this section is to prove a monotonivity formula  which links the area, the diameter and the Willmore functional of a surface $\Sigma \hookrightarrow (\Rtre, \d+h)$. The surface $\Sigma$ can be seen as immersed in two different Riemannian manifolds: $(\Rtre,\d)$ and $(\Rtre, \d+h)$. It follows that all the geometric quantities can be computed with respect the two different spaces and will have different values: the euclidean and the perturbed ones. We use the convention that all the quantities computed with respect to the euclidean metric will have a subscript ''$e$'', for example $|\Sigma|_e, (A_e)_{ij}, H_e, W_e(\Sigma), \ldots$ will denote the euclidean area of $\Sigma$, euclidean second fundamental form, euclidean mean curvature, euclidean Willmore functional, and the corresponding ones evaluated in perturbed metric will have a subscript ''$h$'', for example $|\Sigma|_h, (A_h)_{ij}, H_h, W_h(\Sigma), \ldots$ are the corresponding quantities in metric $\d+h$. Let us start with a straightforward but useful lemma.

\begin{lem} \label{lem:Dist}
Assume that $\|h\|_{C^0(\Rtre)}\leq\eta<1$. It follows that
\begin{itemize}
\item[i)] $(\Rtre,\d+h)$ is a complete Riemannian manifold,
\item[ii)] for every pair of points $p_1,p_2\in \Rtre$ we have
$$ \frac{1}{\sqrt{1+\eta}}d_h(p_1,p_2)\leq |p_1-p_2|_{\Rtre}\leq \frac{1}{\sqrt{1-\eta}} d_h(p_1,p_2),$$
where $|p_1-p_2|_{\Rtre},d_h(p_1,p_2)$ denote the distance respectively in $(\Rtre,\d)$ and in $(\Rtre,\d+h)$ between $p_1$ and $p_2$.
\end{itemize}
\end{lem}

\begin{pf}
To get $i)$ it is sufficient to prove that all the geodesics of $(\Rtre,\d+h)$ are defined on the whole $\R$. Consider the geodesic differential equation $\ddot{x}^\mu+\Gamma^{\mu}_{\nu\lambda} \dot{x}^\nu \dot{x}^{\lambda}=0$ and observe that the Christoffel symbols  $\Gamma^{\mu}_{\nu\lambda}$ of $(\Rtre,\d+h)$ are bounded. Since the geodesics of $(\Rtre,\d+h)$ can be parametrized by arclength, the geodesic differential equation can be interpreted as a dynamical system on the spherical bundle $S(\Rtre,\d+h)$ of $(\Rtre,\d+h)$ (the bundle of the unit tangent vectors) generated by the vector field $X_h(x^\mu,y^\mu):= (y^\mu,-\Gamma^{\mu}_{\nu\lambda} y^\nu y^\lambda)$, where $x\in \Rtre, y \in T_x\Rtre$ with $|y|_h=1$. But $X_h$ is a bounded vector field on $S(\Rtre,\d+h)$ which implies by standard ODE arguments (see for instance Lemma 7.2 and Lemma 7.3 of \cite{AM}) that the integral curves are defined on the whole $\R$. \\
\\
For $ii)$ consider the segment of the straight line $[p_1,p_2]$ connecting $p_1$ and $p_2$. Then by definition we have
$$d_h(p_1,p_2)\leq length_h([p_1,p_2])=\int_0^1\sqrt{(\d+h)(p_2-p_1, p_2-p_1)}\leq \sqrt{1+\eta}|p_1-p_2|_{\Rtre},$$
where of course $length_h([p_1,p_2])$ is the length of the segment $[p_1,p_2]$ in the metric $\d+h$.\\
\\
On the other hand let $\gamma_h:[0,1]\to \Rtre$ be a minimizing geodesic in $(\Rtre, \d+h)$ connecting $p_1$ and $p_2$ (it exists by part $i)$). Then
$$d_h(p_1,p_2)=\int_0^1 \sqrt{(\d+h)(\dot{\gamma_h},\dot{\gamma_h})}\geq \int_0^1 \sqrt{(1-\eta)}\,|\dot{\gamma_h}|_{\Rtre}= \sqrt{1-\eta}\,\,length_e(\gamma_h)\geq \sqrt{1-\eta}\,|p_1-p_2|_{\Rtre},$$
where of course $length_e(\gamma_h)$ is the length of $\gamma_h$ in euclidean metric.
\end{pf}

\begin{lem}\label{lem:AreaEst}
Let $\Sigma \hookrightarrow \Rtre$ be an immersed, smooth, closed, orientable surface, and let $||h||_{C^0}\leq\eta<1/4$. The first fundamental form induced on $\Sigma$ by the two different metrics will be denoted respectively by $\mathring{\d}_{ij}$ and $(\mathring{\d+h})_{ij}$ or simply by $\mathring{\d}$ and $(\mathring{\d+h})$. Then the following pointwise estimate for the area form holds:
\begin{equation} \label{eq:EstAreaForm}
(1-4\eta) \sqrt{\det\mathring{\d}} \leq \sqrt{\det(\mathring{\d+h})} \leq (1+4\eta) \sqrt{\det\mathring{\d}}.
\end{equation} 
%Denoted with $B_\rho^e(x)$ the euclidean ball of radius $\rh$ and center $x$, we will call $\Sigma_{x,\rho}:=\Sigma\cap B_\rho^e(x)$. Then just integrating one gets 
%$$ (1-4\eta) |\Sigma_{x,\rh}|_e \leq |\Sigma_{x,\rh}|_h \leq (1+4\eta) |\Sigma_{x,\rh}|_e $$
%for all $x\in \Rtre$ and $\rh>0$.
\end{lem}

\begin{pf}
Let $f:\Omega\subset \Rdue \to \Rtre$ be a coordinate patch for the surface $\Sigma$. Of course it is enough to do all the computation for a general patch; moreover we can assume that the patch is conformal with respect to the euclidean metric (i.e. we are using isothermal coordinates w.r.t. the euclidean structure).
%Recall that
%$$|\Sigma_{x,\rh}|_e:= \int_{\Sigma_{x,\rh}} \sqrt{\det (\mathring{\d})}$$
%and
%$$|\Sigma_{x,\rh}|_h:= \int_{\Sigma_{x,\rh}} \sqrt{\det (\mathring{\d+h})}$$
%where $\mathring{\d}$ and $\mathring{\d+h}$ are the first fundamental forms induced by  euclidean and perturbed metric.
%\\Let $f_i, i=1,2$ be the derivatives of $f$ with respect to the two coordinates (i.e. the two tangent vectors of the coordinate frame), then by definition:
By definition we have
$$(\mathring{\d+h})_{ij}=(\d+h)(\partial_i f,\partial_j f)=\mathring{\d}_{ij}+h(\partial_i f,\partial_j f).$$
By the choice of the coordinate patch we have that $\mathring{\d}_{ij}$ is diagonal. It follows that
\begin{equation} \label{eq:det}
\det(\mathring{\d+h})=\det(\mathring{\d})+\mathring{\d}_{11} h(\partial_2 f,\partial_2 f)+\mathring{\d}_{22} h(\partial_1 f,\partial_1 f)+\det (h(\partial_i f,\partial_j f)).
\end{equation}
By assumption and by Schwartz inequality we have
$$|h(\partial_i f,\partial_i f)|\leq \eta \mathring{\d}_{ii},\quad h(\partial_1 f,\partial_2 f)^2\leq \eta^2 \mathring{\d}_{11}\mathring{\d}_{22}.$$
Putting these estimates in \eqref{eq:det} and observing that $\eta^2<\eta$ we get
\begin{equation} \label{eq:EstDet}
(1-4\eta) (\det \mathring{\d}) \leq \det (\mathring{\d+h}) \leq (1+4\eta) (\det \mathring{\d}),
\end{equation}
and the lemma follows.
%Since in our range $1-4\eta \leq \sqrt{1-4\eta}$ and $\sqrt{1+4\eta} \leq 1+4\eta$, we have the thesis just taking the square root of \eqref{eq:EstDet} and integrating on the desired domain.
\end{pf}

In the following lemma we derive a pointwise estimate from above and below of the mean curvature squared in perturbed setting in terms of the corresponding euclidean quantities.

\begin{lem}\label{lem:estH}
Let $\Sigma \hookrightarrow \Rtre$ be an immersed, smooth, closed, orientable surface. Assume that $\|h\|_{C^0}\leq \eta$ and $\|Dh\|_{C^0}\leq \theta$ with $\eta$ small. Then the following pointwise estimate holds:
$$(1-C\eta-\gamma) |H_e|^2-(C\eta+\gamma) |A_e|^2-C_\gamma \th^2 \leq |H_h|^2 \leq (1+C\eta+\gamma) |H_e|^2+(C\eta+\gamma) |A_e|^2+C_\gamma \th^2,$$
where $\gamma >0$ is arbitrary and $C_\gamma\leq C(1+\frac{1}{\gamma})$.
\end{lem}

\begin{pf}
Let $p\in \Sigma$ and choose the parametrization $f$ given by the normal coordinates at $p$ with respect to the metric $\mathring{\delta}$, such that the coordinate vectors $\partial_i f$ are euclidean-orthonormal and diagonalize the euclidean second fundamental form $A_e$ at $p$ (the first condition is trivial, the second can be achieved by a rotation). With this choice of coordinates, the euclidean Christoffel symbols $\tilde{\Gamma}^k_{ij}$ of $\Sigma$ vanish at $p$ and therefore
\begin{equation}\label{eq:D2}
\partial^2_{ij} f(p)=(A_e)_{ij}(p) \nu_e(p) + \tilde{\Gamma}^k_{ij}(p)\partial_k f(p)= (A_e)_{ij}(p) \nu_e(p),
\end{equation}
where $\nu_e$ denotes the euclidean normal vector to $\Sigma$, namely $\nu_e= \partial_1 f\times\partial_2 f$.\\
\\
The normal vector to $\Sigma$ in perturbed metric is denoted $\nu_h$ and has the form $\nu_h=\nu_e+N$, where the correction $N$ is small since $\|h\|_{C^0}$ is small. More precisely it follows from the orthogonality conditions $(\d+h)(\partial_i f,\nu_h)=0$ that
$$\d(\partial_i f,N)=-h(\partial_i f,\nu_e)+\text{higher order terms}.$$
Imposing the normalization condition $(\delta+h)(\nu_h,\nu_h)=1$ we obtain
\begin{equation}
\d(N,\nu _e)= - \frac{1}{2}h(\nu _e ,\nu _e)+\text{higher order terms}.\nonumber
\end{equation}
Since $(\partial_1 f,\partial_2 f,\nu_e)$ is an orthonormal frame in euclidean metric, we can represent $N$ as 
\begin{equation}\label{eqN}
N= -h(\nu_e ,\partial_1 f)\partial_1 f-h(\nu_e,\partial_2 f)\partial_2 f-\frac{1}{2}h(\nu_e,\nu_e)\nu_e+\text{higher order terms}.
\end{equation}
Observe that the higher order terms can be computed in an inductive way using the orthonormalization conditions above and that for $\eta$ small
\begin{equation}\label{eq:N}
|N|_e=\sqrt{\d(N,N)}\leq C\eta. 
\end{equation} 
Now let us compute the perturbed second fundamental form 
$$(A_h)_{ij}=(\d+h)(\nu_h,\,^{\d+h}\nabla_{\partial_i f}\partial_j f),$$
where $^{\d+h}\nabla $ is the covariant derivative in $(\Rtre,\d+h)$. By definition 
$$^{\d+h}\nabla_{\partial_i f}\partial_j f=\partial^2_{ij} f+\,^{\d+h}\Gamma\partial_i f\partial_j f,$$ 
where $ ^{\d+h}\Gamma$ are the Christoffel symbols of $(\Rtre, \d+h)$
and $^{\d+h}\Gamma\partial_i f\partial_j f=\,^{\d+h}\Gamma^\mu_{\nu \lambda} \partial_i f^\nu \partial_j f^\lambda e_\mu$, where $\{e_\mu\}$ is the standard euclidean orthonormal basis of $(\Rtre,\d)$ and $\partial_i f=\partial_i f^\mu e_\mu$.
\\Using \eqref{eq:D2}, the perturbed second fundamental form becomes
$$(A_h)_{ij}=(\d+h)(\nu_e+N, (A_e)_{ij} \nu_e+ \,^{\d+h}\Gamma \partial_i f\partial_j f). $$
Observing that $|^{\d+h}\Gamma|\leq C\theta$ and recalling \eqref{eq:N} one gets
\begin{equation} \label{eq:EstA}
(A_e)_{ij}-C\eta(A_e)_{ij}-C\theta \leq(A_h)_{ij}\leq (A_e)_{ij}+C\eta(A_e)_{ij}+C\theta.
\end{equation}
Squaring and using the $\gamma$-Cauchy inequality we get that for any $\gamma>0$ 
\begin{equation}\label{eq:EstAe2Ah2}
(1-2\gamma-C\eta)|A_e|^2-C_{\gamma}\theta^2(1+\eta^2) \leq |A_h|^2\leq (1+2\gamma+C\eta)|A_e|^2+C_{\gamma}\theta^2(1+\eta^2),
\end{equation}
where $C_\gamma\leq C(1+\frac{1}{\gamma})$. Since $H_h=(\mathring{\d+h})^{ij} (A_h)_{ij}$ by definition and since
$$(\mathring{\d})^{ij}-C\eta \leq (\mathring{\d+h})^{ij} \leq (\mathring{\d})^{ij}+C\eta,$$
we get by taking the trace in \eqref{eq:EstA} with respect to $\mathring{\d+h}$ that
\begin{equation}\label{eq:EstH}
H_e-C\eta |A_e|_e-C\theta \leq H_h \leq H_e+C\eta |A_e|_e+C\theta,
\end{equation}
where $|A_e|_e$ (in the sequel called just $|A_e|$) is the euclidean norm of the euclidean second fundamental form. Using the Cauchy inequality it follows that
\begin{eqnarray}
|H_h|^2&\leq&|H_e|^2+C\eta |H_e| |A_e|+C\theta |H_e|+C\eta^2 |A_e|^2+C\theta |A_e|+C\theta^2\nonumber\\
&\leq&(1+C\eta+\gamma) |H_e|^2+ (C\eta+\gamma) |A_e|^2+C_\gamma \theta^2. \nonumber
\end{eqnarray}
The estimate from below is analogous, and the lemma is proved. 
\end{pf}

%\begin{lem}\label{lem:EstW}
%Let  $\spt h\subseteq B^e_{r_0}(x_0)$ for some $x_0\in \Rtre$ and $r_0>0$. As before $\|h\|_{C^0}\leq \eta$, $\|Dh\|_{C^0}\leq \theta$ ($\eta$ is supposed to be small while no assumption is made on $\theta$) and $\Sigma\hookrightarrow \Rtre$ is a closed smooth orientable surface of genus $0$ (it is enough to ask the uniform bound  $\gen (\Sigma_k)\leq g$) immersed in $\Rtre$. Then
%\begin{equation}\label{eq:EstW}
%(1-C\eta-C\gamma-C_\gamma r_0^2 \theta^2) W_e(\Sigma)- C_g(\eta+\gamma)\leq W_h(\Sigma)
%\end{equation}
%where $C_g\to \infty$ as $g\to \infty$, $\gamma >0$ can be chosen arbitrarily small and $C_\gamma$ is a constant depending on $\gamma$ such that $C_\gamma\to \infty$ if $\gamma\to 0$ but which can be bounded by $C_\gamma\leq C(1+\frac{1}{\gamma})$ for $C$ large enough independent on $\gamma$. It follows that for $\gamma,\eta$ and $r_0\theta$ small enough
%\begin{equation}\label{eq:EstW1}
%W_e(\Sigma) \leq \frac{3}{2} W_h(\Sigma)+1.
%\end{equation}
%\end{lem}

\begin{lem}\label{lem:EstW}
Let $\Sigma\hookrightarrow \Rtre$ be an immersed, smooth, closed, orientable surface. Assume that $\|h\|_{C^0}\leq \eta$ ($\eta>0$ small) and $\|Dh\|_{C^0}\leq \theta$, and that  $\spt h\subset B^e_{r_0}(x_0)$  where  $B^e_{r_0}(x_0)$ is the euclidean ball of center $x_0\in \Rtre$ and radius $r_0>0$. Then
\begin{equation}\label{eq:EstW}
(1-C\eta-C\gamma-C_\gamma r_0^2 \theta^2) W_e(\Sigma)- C_g(\eta+\gamma)\leq W_h(\Sigma),
\end{equation}
where $\gamma>0$ is arbitrary, $C_g\leq C(1+\gen \Sigma)$ is a constant depending on $\gen\Sigma$ and $C_\gamma\leq C(1+\frac{1}{\gamma})$. Moreover it follows for $\eta$ and $r_0\theta$ sufficiently small that
\begin{equation}\label{eq:EstW1}
W_e(\Sigma) \leq \frac{3}{2} W_h(\Sigma)+1.
\end{equation}
\end{lem}

\begin{pf}
Recalling the estimate of the area form \eqref{eq:EstAreaForm}, integrating the formula of Lemma \ref{lem:estH} yields
$$W_h(\Sigma)=\frac{1}{4}\int_{\Sigma}|H_h|^2\sqrt{\det \mathring{(\d+h)}} \geq \int_{\Sigma} \Big[\Big(\frac{1}{4}-C\eta-\gamma \Big) |H_e|^2-(C\eta+\gamma) |A_e|^2-C_\gamma \th^2\chi_h\Big] (1-4\eta) \sqrt{\det\mathring{\d}},$$
where $\chi_h$ is the characteristic function of $\spt h$. From the Gauss-Bonnet Theorem it follows that
$$\int_\Sigma |A_e|^2\sqrt{\det\mathring{\d}}=\int_\Sigma |H_e|^2\sqrt{\det\mathring{\d}}-4\pi \chi_E (\Sigma),$$
where $\chi_E(\Sigma)=2-2\gen\Sigma$ is the Euler characteristic of $\Sigma$. %In the case $\gen (\Sigma)=0$ of course $\chi_E(\Sigma)=2$ but more generally if $\gen (\Sigma)$ is uniformly bounded also $-4\pi\chi_E(\Sigma)$ will be uniformly bounded from above. 
Hence
$$W_h(\Sigma)\geq (1-C\eta-C\gamma) W_e(\Sigma)-C_g(\eta+\gamma)-C_\gamma \th^2 |\Sigma\cap \spt h|_e.$$
From formula (1.3) in \cite{SiL} it follows that
\begin{equation}\label{eq:AreaW}
|\Sigma\cap \spt h|_e\leq |\Sigma\cap B^e_{r_0}(x_0)|_e\leq C r_0^2 W_e(\Sigma).
\end{equation}
%We can conclude that
%$$W_h(\Sigma)\geq (1-C\eta-\gamma-C_\gamma\th^2 r_0^2) W_e(\Sigma)-C_g(\eta+\gamma).$$
%We get the thesis by first fixing $\gamma$ small enough and then choosing sufficiently small $\eta,\th$.
Therefore the lemma is proved.
\end{pf}

Using the estimates of the previous lemmas, we get the desired monotonicity formula in the following proposition. For that we define $\Sigma_{x,\rho}:=\Sigma\cap B_\rho^e(x)$.

%\begin{pro}\label{prop:MFh}
%As before let  $\spt h\subseteq B^e_{r_0}(x_0)$ for some $x_0\in \Rtre$, $r_0>0$ and $\|h\|_{C^0}\leq \eta$, $\|Dh\|_{C^0}\leq \theta$; recall that $\Sigma\hookrightarrow \Rtre$ is a closed smooth orientable surface of genus $g\geq 0$ immersed in $\Rtre$. Then for $\eta$ and $r_0\th$ small enough the following  inequality holds
%$$\sigma^{-2}|\Sigma_{x,\sigma}|_h\leq C \Big[\rh^{-2}|\Sigma_{x,\rh}|_h+ W_h(\Sigma_{x,\rh})+ [C_g(\eta+\gamma)+C_\gamma r_0^2 \th^2](W_h(\Sigma)+1)  \Big] \quad 0<\sigma\leq\rh<\infty $$
%where $\gamma >0$ can be chosen arbitrarily small and $C_\gamma, C_g$ are constants depending on $\gamma$ (respectively on $g$) such that $C_\gamma\to \infty$ if $\gamma\to 0$ (respectively $C_g\to \infty$ if $g\to \infty$). It follows the more simple estimate
%\begin{equation}\label{eq:MF}
%\sigma^{-2}|\Sigma_{x,\sigma}|_h\leq C_g \Big[\rh^{-2}|\Sigma_{x,\rh}|_h+ W_h(\Sigma)+1  \Big]\quad 0<\sigma\leq\rh<\infty
%\end{equation}
%and just taking the limit $\rh\to \infty$,
%\begin{equation}\label{eq:EstAreah}
%|\Sigma|_h\leq C_g(W_h(\Sigma)+1) (\diam_e\Sigma)^2
%\end{equation}
%where $\diam_e\Sigma$ is the euclidean diameter of $\Sigma$.
%\end{pro}

\begin{pro}\label{prop:MFh}
Let $\Sigma\hookrightarrow \Rtre$ be an immersed, smooth, closed, orientable surface. Assume $\|h\|_{C^0}\leq \eta$ and $\|Dh\|_{C^0}\leq \theta$, and that  $\spt h\subset B^e_{r_0}(x_0)$ for some $x_0\in \Rtre$ and $r_0>0$. Then for $\gamma$, $\eta$ and $r_0\theta$ sufficiently small the following inequality holds
$$\sigma^{-2}|\Sigma_{x,\sigma}|_h\leq C \Big[\rh^{-2}|\Sigma_{x,\rh}|_h+ W_h(\Sigma_{x,\rh})+ [C_g(\eta+\gamma)+C_\gamma r_0^2 \th^2](W_h(\Sigma)+1)  \Big] \quad\text{for all }0<\sigma\leq\rh<\infty, $$
where $C_g\leq C(1+\gen \Sigma)$ is a constant depending on $\gen\Sigma$ and $C_\gamma\leq C(1+\frac{1}{\gamma})$. 
%Moreover it follows
%\begin{equation}\label{eq:MF}
%\sigma^{-2}|\Sigma_{x,\sigma}|_h\leq C_g \Big[\rh^{-2}|\Sigma_{x,\rh}|_h+ W_h(\Sigma)+1  \Big]\quad\text{for all }0<\sigma\leq\rh<\infty,
%\end{equation}
%Moreover by letting $\rh\to \infty$ we get
%\begin{equation}\label{eq:EstAreah}
%|\Sigma|_h\leq C_gC_\gamma(W_h(\Sigma)+1) (\diam_e\Sigma)^2,
%\end{equation}
%where $\diam_e\Sigma$ is the euclidean diameter of $\Sigma$.
\end{pro}

\begin{pf}
Let us recall the euclidean monotonicity formula proved by Simon (formula (1.3) in \cite{SiL}):
\begin{equation}\label{EMF}
\sigma^{-2} |\Sigma_{x,\sigma}|_e\leq C(\rh^{-2}|\Sigma_{x,\rh}|_e+W_e(\Sigma_{x,\rh})) .
\end{equation}
We just have to estimate from above and below the area part and from above the Willmore term. From Lemma \ref{lem:AreaEst} it follows by integration that
$$\frac{1}{1+4\eta} |\Sigma_{x,\sigma}|_h\leq |\Sigma_{x,\sigma}|_e,\quad|\Sigma_{x,\rh}|_e\leq \frac{1}{1-4\eta} |\Sigma_{x,\rh}|_h.$$
Integrating the formula of Lemma \ref{lem:estH} yields
$$W_h(\Sigma_{x,\rh}) \geq \int_{\Sigma_{x,\rh}} \Big[\Big(\frac{1}{4}-C\eta-\gamma\Big) |H_e|^2-(C\eta+\gamma) |A_e|^2-C_\gamma \th^2\chi_h\Big] (1-4\eta) \sqrt{\det\mathring{\d}},$$
where again $\chi_h$ is the characteristic function of $\spt h$. From the Gauss Bonnet Theorem and \eqref{eq:EstW1} we get
$$\int_{\Sigma_{x,\rh}}|A_e|^2\sqrt{\det\mathring{\d}}\leq \int_{\Sigma}|A_e|^2\sqrt{\det\mathring{\d}}\leq C_g(W_e(\Sigma)+1)\leq C_g(W_h(\Sigma)+1),$$
where $C_g\leq C(1+\gen \Sigma)$ is a constant depending on $\gen\Sigma$. Hence
$$W_h(\Sigma_{x,\rh}) \geq (1-C\eta-C\gamma) W_e(\Sigma_{x,\rh})- C_g(\eta+\gamma)(W_h(\Sigma)+1) -C_\gamma \th^2|\Sigma_{x,\rh}\cap \spt h|_e.$$
As before 
$$|\Sigma_{x,\rh}\cap \spt h|_e\leq |\Sigma\cap B^e_{r_0}(x_0)|_e\le C r_0^2 W_e(\Sigma)\leq C r_0^2 (W_h(\Sigma)+1),$$
and thus we get for $\eta$ and $\gamma$ sufficiently small that
$$W_e(\Sigma_{x,\rh}) \leq C W_h(\Sigma_{x,\rh})+ C_g(\eta+\gamma)(W_h(\Sigma)+1) +C_\gamma r_0^2\th^2 (W_h(\Sigma)+1),$$
%From the previous inequalities \eqref{eq:EstAreaForm}-\eqref{eq:AreaW} and \eqref{eq:EstW1} 
%$$|\Sigma_{x,\rh}\cap \spt h|_h\leq C |\Sigma\cap B^e_{r_0}(x_0)|_e\leq  C r_0^2 W_e(\Sigma)\leq C r_0^2 (W_h(\Sigma)+1),$$
%hence
%$$W_e(\Sigma_{x,\rh}) \leq C W_h(\Sigma_{x,\rh})+ C_g(\eta+\gamma)(W_h(\Sigma)+1) +C_\gamma r_0^2 \th^2 (W_h(\Sigma)+1)$$
%and we can conclude that
%$$\sigma^{-2}|\Sigma_{x,\sigma}|_h\leq C \Big[\rh^{-2}|\Sigma_{x,\rh}|_h+ W_h(\Sigma_{x,\rh})+ [C_g(\eta+\gamma)+C_\gamma r_0^2 \th^2](W_h(\Sigma)+1)  \Big].$$
and the proposition follows from Simon's monotonicity formula \eqref{EMF}.
\end{pf}

\subsection{A priori estimates for a minimizing sequence of $W$}\label{Sec:AprioriEst}

Under a very general assumption on the metric (we ask that the scalar curvature of the ambient manifold is strictly positive in one point) we will show global a priori estimates for minimizing sequences of the Willmore functional; more precisely we get uniform upper area bounds, uniform upper and lower bounds on the diameters and we show that minimizing sequences are contained in a compact subset of $\Rtre$.

\begin{pro}\label{pro:a}
Following the previous notation, assume that the scalar curvature $R_h$ of $(\Rtre,\delta+h)$ is strictly positive in some point $\overline{x} \in \Rtre$, namely $R_h(\overline{x}) > 0$, then
$$\inf_{f\in[\Sp^2,(\Rtre,\delta+h)]}W_h(f)<4\pi.$$
\end{pro}

\begin{pf}
From Proposition 3.1 of \cite{Mon1}, on geodesic spheres $S_{\overline{x},\rh}$ of center $\overline{x}$ and small radius $\rh$ one has 
$$W_h(S_{\overline{x},\rh})=4\pi-\frac{2\pi}{3}R_h(\overline{x})\rh^2+O(\rh^3).$$
Since these surfaces are smooth embeddings of $\Sp^2$ and $R_h(\overline{x})>0$, the conclusion follows. 
\end{pf}

%\begin{cor}\label{cor:sk}
%Let $\Sigma_k$, with  $\gen \Sigma_k\leq g$, be a minimizing sequence for $W_h$ of bounded genus
%$$W_h(\Sigma_k)\downarrow \a_h^g$$ 
%and assume there exists a point $\bar{p}\in\Rtre$ such that the scalar curvature $R_h(\bar{p})>0$.
%\\Then there exists $\e>0$ such that for large $k$ 
%$$W_h(\Sigma_k)<4\pi-\e.$$
%\end{cor}

%Now let us state and prove uniform a priori upper bounds on the minimizing sequence $\Sigma_k$. The idea is to use just that $W_h(\Sigma_k)<4\pi-\e$ and then perform a blow down procedure making use of the rescaling invariance of the Willmore functional (see equation \eqref{eq:ScaleInv} below).

The last proposition together with (\ref{eq:EstW1}) implies  that if the scalar curvature $R_h$ of $(\Rtre,\delta+h)$ is strictly positive in some point, then for a minimizing sequence $f_k\in[\Sp^2,(\Rtre,\delta+h)]$ of the functional $W_h$ in $[\Sp^2,(\Rtre,\delta+h)]$ we have for $k$ sufficiently large
$$W_e(f_k)<8\pi,$$
and thus $f_k$ is an embedding. Therefore in order to minimize the functional $W_h$ in $[\Sp^2,(\Rtre,\delta+h)]$ we can take minimizing sequences of smooth spheres $\Sigma_k$ embedded in $\Rtre$.

%\begin{pro}\label{pro:Placebound}
%Let $(\Rtre, \d+h)$ be as before with small $\|h\|_{C^1}$  and let $\Sigma_k \hookrightarrow \Rtre$ be a sequence of immersed smooth closed orientable surfaces of genus 0 (more generally one can ask the uniform bound  $\gen (\Sigma_k)\leq g$, but in this case the required smallness of $\|h\|_{C^1}$ depends on $g$ and goes to $0$ as $g\to \infty$). Assume that
%$$\limsup _k W_h(\Sigma_k)<4\pi,$$
%then  
%\\i)there exists a compact $K \subset \Rtre$ such that 
%$$ \Sigma_k\subseteq K \quad \forall k\in \N,$$
%ii) there exists a uniform area bound  
%$$|\Sigma_k|_h\leq C. $$
%for some large $C>0$.
%\end{pro}

\begin{pro}\label{pro:Placebound}
Assume $\|h\|_{C^0}\leq \eta$ and $\|Dh\|_{C^0}\leq \theta$, and  that  $\spt h\subset B^e_{r_0}(x_0)$ for some $x_0\in \Rtre$ and $r_0>0$.  
Suppose  that $\inf_{f\in[\Sp^2,(\Rtre,\delta+h)]}W_h(f)<4\pi$ and let $\Sigma_k \hookrightarrow \Rtre$ be a minimizing sequence of smooth, embedded spheres for the functional $W_h$ in $[\Sp^2,(\Rtre,\delta+h)]$. Then for  $\eta$ and $r_0\theta$ sufficiently small we have that
\begin{itemize}  
\item[i)] there exists a compact set $K \subset \Rtre$ such that $\Sigma_k\subset K$ for $k$ sufficiently large,
\item[ii)] there exists a constant $C<\infty$ such that $|\Sigma_k|_h\leq C$ for $k$ sufficiently large.
\end{itemize}
\end{pro}

\begin{pf}
First of all observe that each surface $\Sigma_k$ is connected. As before let $\eta= \|h\|_{C^0}$ and $\th=\|Dh\|_{C^0}$, and let $r_0>0$ such that $\spt h \subset B_{r_0}^e(0)$. From $W_h(\Sigma_k)<4\pi$ it follows that 
$$\Sigma_k \cap B_{r_0}^e(0)\neq \emptyset,$$
since otherwise $W_h(\Sigma_k)=W_e(\Sigma_k)$ and thus $W_e(\Sigma_k)\geq 4\pi$ by Theorem 7.2.2 in \cite{Will}.\\
\\
The goal is to prove that $\limsup_k(\diam_e \Sigma_k)<\infty$, because then $i)$ follows immediately, and statement $ii)$ follows by letting $\rho\to\infty$ in Proposition \ref{prop:MFh}. Assume that up to subsequences
$$\diam_e \Sigma_k \nearrow \infty.$$
For each $k$ we rescale in the following way. We set
\begin{equation}\label{def:stilde}
\tilde{\Sigma}_k=\frac{1}{\diam_e \Sigma_k}\Sigma_k,\quad(h_k)_{\mu\nu}(x)=h_{\mu\nu}\big((\diam_e \Sigma_k)\,x\big)%\label{def:hk}
\end{equation}
It follows that
\begin{equation}\label{eq:diamk}
\diam_e \tilde{\Sigma}_k=1,\quad\spt h_k=\frac{1}{\diam_e \Sigma_k} \spt h \subseteq B_{r_k}^e(0), %\label{eq:suppk}
\end{equation}
where 
\begin{equation}\label{def:rk}
r_k=\frac{1}{\diam_e \Sigma_k} r_0 \searrow 0.
\end{equation}
Let $\eta_k=\|h_k\|_{C^0}$ and $\th_k=\|Dh_k\|_{C^0}$, and observe that 
\begin{equation}\label{eq:etak}
\eta_k=\eta,\quad r_k \th_k=r_0 \th. %\label{eq:rth}
\end{equation}
%The second equality follows simply from the chain rule
%$$\frac{\partial}{\partial x^\lambda} (h_k)_{\mu\nu}|_{x}= \frac{\partial}{\partial x^\lambda} [h_{\mu\nu}(\diam_e \Sigma_k \; .)]|_{x}=\diam_e \Sigma_k \; \frac{\partial}{\partial x^\lambda} h_{\mu\nu}|_{\diam_e \Sigma_k \; x}.$$
Moreover, just from the definitions, it is easy to check the scale invariance of the Willmore functional
\begin{equation}\label{eq:ScaleInv}
W_{h_k}(\tilde{\Sigma}_k)=W_h(\Sigma_k).
\end{equation}
Because of (\ref{eq:etak}), for $\eta$ and $r_0 \theta$ sufficiently small we can apply   Proposition \ref{prop:MFh} to $\tilde{\Sigma}_k$ to get in view of (\ref{eq:diamk}) and the uniform bound on the Willmore energy of $\tilde{\Sigma}_k$ that
\begin{equation}\label{eq:Areaboundk}
|\tilde{\Sigma}_k|_{h_k}\leq C.
\end{equation}  
Now it follows from \eqref{eq:EstW1} and Lemma \ref{lem:AreaEst} that
\begin{equation}\label{eq:BoundAe}
|\tilde{\Sigma}_k|_e\leq C,\quad W_{e}(\tilde{\Sigma}_k)\leq C.
\end{equation} 
Now define the integral, rectifiable 2-varifold $\mu_k^e$ in $(\Rtre, \d)$ by
\begin{equation}\label{3}
\mu_k^e=\mathcal{H}^2_e\llcorner\tilde\Sigma_k,
\end{equation}
where $\mathcal{H}^2_e$ denotes the usual 2-dimensional Hausdorff measure. It follows that $\mu_k^e(\Rtre)\le C$ and that the first variation can be bounded by a universal constant by (\ref{eq:BoundAe}). By a compactness result for varifolds (see \cite{SiGMT}), there exists an integral, rectifiable 2-varifold $\mu^e$ in $(\Rtre, \d)$ with weak mean curvature vector $H_e\in L^2(\mu^e)$, such that (after passing to a subsequence) $\mu_k^e\to\mu^e$ weakly as measures and 
\begin{equation}\label{4'}
W_e(\mu^e)=\frac{1}{4}\int|H_e|^2\,d\mu^e\le\liminf_{k\to\infty}W_e(\tilde{\Sigma}_k)\le C.
\end{equation}
More precisely we have the following: For fixed $n$ we have due to (\ref{eq:diamk}) that $\spt h_k\subset B_\frac{1}{n}(0)$ for $k$ sufficiently large. It follows from the varifold convergence, the lower semicontinuity of the Willmore functional, the assumption and (\ref{eq:ScaleInv}) that
\begin{equation}
W_e(\mu^e\llcorner\Rtre\setminus B_\frac{1}{n}(0))\le\liminf_{k\to\infty}W_e(\tilde{\Sigma}_k\setminus B_\frac{1}{n}(0))=\liminf_{k\to\infty}W_{h_k}(\tilde{\Sigma}_k\setminus B_\frac{1}{n}(0))\le\liminf_{k\to\infty}W_h(\Sigma_k)<4\pi.
\end{equation}
Since $H_e\in L^2(\mu^e)$ it follows by letting $n\to\infty$ that
\begin{equation}\label{<4pi}
W_e(\mu^e)<4\pi.
\end{equation}
Now we want to prove that actually $\mu^e$ is not the null varifold. For that we will prove that there exists a $\beta>0$ such that 
\begin{equation}\label{eq:TesiClaim}
\mu_k^e(B_1^e(0)\backslash B^e_\frac{1}{2}(0)) \ge \beta \quad\text{for large $k$},
\end{equation}
because then it would follow from the weak convergence that
%Since the varifold convergence implies the weak convergence of measures of the associated mass measures, we will have for the converging subsequence
\begin{equation}\label{V2beta}
\mu^e\Big(\overline{B_1^e(0)\backslash B^e_\frac{1}{2}(0)}\Big)\ge \limsup_{k\to\infty}\mu_k^e(B_1^e(0)\backslash B^e_\frac{1}{2}(0))\ge\beta.
\end{equation}
%for some $\varepsilon\geq 0$ such that $\|V\|(\partial B^e_{1/2-\varepsilon}(0))=0$ (it exists since $\|V\|$ is a finite measure); observe we denoted with $\|V\|(\Rtre\backslash B^e_{1/2}(0))$  the measure of $\Rtre\backslash B^e_{1/2}(0)$ with respect to the mass measure $\|V\|$ of the varifold $V$. This will prove the Claim.
To prove \eqref{eq:TesiClaim}, notice that, since $\tilde{\Sigma}_k$ is connected, $\diam_e \tilde{\Sigma}_k = 1$, $\tilde{\Sigma}_k\cap B^e_{r_k}(0)\neq \emptyset$ and $r_k\to0$, it follows that 
\begin{equation}\label{eq:supp}
 \spt h_k \subset B^e_\frac{1}{2}(0),\quad\tilde{\Sigma}_k \cap \partial B^e_\frac{3}{4}(0) \neq \emptyset \quad  \text{for $k$ sufficiently large.}
\end{equation} 
For $N\in\N$ let
$$A_i=B^e_{\frac{1}{2}+\frac{i}{4N}}(0)\backslash B^e_{\frac{1}{2}+\frac{(i-1)}{4N}}(0), \quad i=1,\ldots, N,$$ 
and observe that $A_i\cap A_j=\emptyset$ for $i\neq j$ and that 
$$B^e_\frac{3}{4}(0)\backslash B^e_\frac{1}{2}(0)=\bigcup_{i=1}^N A_i.$$
%Let us do a partition of $B^e_{3/4}(0)\backslash B^e_{1/2}(0)$ with $N$ annuli at distance $\frac{1}{4N}$ one each other, i.e.  the $i^{th}$ annulus is of the type $$A_i:=B^e_{1/2+\frac{i}{4N}}(0)\backslash B^e_{1/2+\frac{(i-1)}{4N}}(0) \quad i=1,\ldots, N.$$ 
Since $\tilde\Sigma_k$ is connected, $\tilde{\Sigma}_k\cap B^e_\frac{1}{2}(0)\supset\tilde{\Sigma}_k\cap B^e_{r_k}(0)\neq \emptyset$ and $\tilde{\Sigma}_k \cap \partial B^e_\frac{3}{4}(0) \neq \emptyset$, it follows that for each $i\in\{1,\ldots,N\}$ there exists a point $x_i^k\in \tilde{\Sigma}_k \cap A_i$ such that $B^e_\frac{1}{8N}(x_i^k)\subset A_i$. Simon's monotonicity formula (formula (1.4) in \cite{SiL}) yields
\begin{equation}\label{mon}
\pi\leq C\big(64 N^2 |\tilde{\Sigma}_k\cap B^e_\frac{1}{8N}(x_i^k) |_e + W_e(\tilde{\Sigma}_k\cap B^e_\frac{1}{8N}(x_i^k)\big).
\end{equation}
Now assume that
$$W_e(\tilde{\Sigma}_k\cap B^e_\frac{1}{8N}(x_i^k))\geq\frac{\pi}{2C}\quad\text{for all }i=1,\ldots,N. $$ 
Since the balls $B^e_\frac{1}{8N}(x_i^k)$, $i=1,\ldots,N$, are pairwise disjoint, we get 
$$W_{e}(\tilde{\Sigma}_k\backslash B_\frac{1}{2}^e(0))\geq \sum _{i=1}^{N} W_e(\tilde{\Sigma}_k\cap B^e_\frac{1}{8N}(x_i^k)) \geq N\frac{\pi}{2C}.$$ 
Since the Willmore energy is uniformly bounded, we get for $N$ sufficiently large a contradiction. Thus there exists a point $x_i^k$ such that $W_e(\tilde{\Sigma}_k\cap B^e_\frac{1}{8N}(x_i^k))\le\frac{\pi}{2C}$, and it follows from (\ref{mon}) that
$$|\tilde{\Sigma}_k\cap B^e_\frac{1}{8N}(x_i^k) |_e \geq \frac{1}{64 N^2} \Big(\frac{\pi}{C}-W_e(\tilde{\Sigma}_k\cap B^e_\frac{1}{8N}(x_i^k) \Big)\ge\frac{\pi}{128 C N^2}>0.$$
This shows (\ref{eq:TesiClaim}).
%We can assume that each $\tilde{\Sigma}_k\cap(\Rtre\backslash B^e_{1/2}(0))$ is connected, otherwise just take a connected component of $\tilde{\Sigma}_k\cap(\Rtre\backslash B^e_{1/2}(0))$ which intersects $S_{3/4}$.
%\\From the connection property, for each annulus and for each $\tilde{\Sigma}_k$ there is a point $x_i^k\in \tilde{\Sigma}_k \cap A_i$ such that $B^e_{1/(8N)}(x_i^k)\subset A_i$.
%From Simon's monotonicity formula (formula (1.4) page 285 of \cite{SiL}),
%$$\pi\leq C\big(64 N^2 |\tilde{\Sigma}_k\cap B^e_{1/(8N)}(x_i^k) |_e + W_e(\tilde{\Sigma}_k\cap B^e_{1/(8N)}(x_i^k)\big).$$
%It follows that
%\begin{equation}\label{eq:lowerarea}
%|\tilde{\Sigma}_k\cap B^e_{1/(8N)}(x_i^k) |_e \geq \frac{1}{64 N^2} \Big( \frac{\pi}{C}-W_e(\tilde{\Sigma}_k\cap B^e_{1/(8N)}(x_i^k) \Big).
%\end{equation}
%Now it is enough to prove  that  $ \exists N$ large enough:  $ \forall k$ large $\exists x_i^k$ (notation above) such that
%$$W_e(\tilde{\Sigma}_k\cap B^e_{1/(8N)}(x_i^k))<\frac{\pi}{2C}.$$
%If it is not true, $\forall N>0$ there exists a large $k$  such that $\forall x_i^k \; i=1,\ldots,N,$
%$$W_e(\tilde{\Sigma}_k\cap B^e_{1/(8N)}(x_i^k))\geq\frac{\pi}{2C}. $$   
%But, for $k$ fixed, the balls $B^e_{1/(8N)}(x_i^k)$ $i=1,\ldots,N$ are disjoint; hence
%$$W_{e}(\tilde{\Sigma}_k\cap(\Rtre\backslash B_{1/2}))\geq \sum _{i=1}^{N} W_e(\tilde{\Sigma}_k\cap B^e_{1/(8N)}(x_i^k)) \geq N\frac{\pi}{2C}.$$ 
%Since $N$ is arbitrarily large, this contradicts the boundness of $W_e(\tilde{\Sigma}_k)\geq W_{e}(\tilde{\Sigma}_k\cap(\Rtre\backslash B_{1/2}))$. 
%\\This concludes the proof of the claim.
%\\\\
Now since $\mu^e\neq0$ is integral, it follows from a generalized monotonicity formula proved by Kuwert and Sch\"atzle in \cite{KS} that $W_e(\mu^e)\geq 4 \pi$, which contradicts (\ref{<4pi}), and thus the proposition is proved.
\end{pf}

Finally we would like to mention that a minimizing sequence $\Sigma_k$ cannot shrink to a point if the scalar curvature $R_h$ of $(\Rtre,\delta+h)$ is strictly positive in some point, namely
\begin{equation}\label{lowdiambound}
\liminf_{k\to\infty}(\diam_h \Sigma_k)>0.
\end{equation} 
This follows from the fact that in this case the infimum of the Willmore energy on the class $[\Sp^2,(\Rtre,\delta+h)]$ is strictly less than $4\pi$ together with Proposition 2.5 in \cite{KMS}, which also holds for non-compact Riemannian manifolds $M$ without boundary, assuming that the minimizing sequence stays in a compact set.

\subsection{Existence and regularity of minimizers for the Willmore energy}\label{Sec:ExistMinim}
Since this semi perturbative setting is closely related to the setting in \cite{SiL}, we just sketch the procedure for proving existence and regularity, pointing out the main differences with \cite{SiL}. We refer to the mentioned paper for more  details and also to \cite{KMS} or \cite{Schy}.\\
\\
Let $\Sigma_k\in[\Sp^2,(\Rtre,\delta+h)]$ be a minimizing sequence for the Willmore energy $W_h$ in perturbed metric. We assume that the scalar curvature $R_h$ of $(\Rtre,\delta+h)$ is strictly positive in some point $\overline{x} \in \Rtre$, namely $R_h(\overline{x}) > 0$. Define the integral, rectifiable 2-varifold $\mu_k^h$ in $(\Rtre, \d+h)$ by
\begin{equation}
\mu_k^h=\mathcal{H}^2_h\llcorner\Sigma_k,
\end{equation}
where $\mathcal{H}^2_h$ is the 2-dimensional Hausdorff measure with respect to the metric $\d+h$. It follows from Proposition \ref{pro:Placebound} and the minimizing sequence property that for $\eta$ and $r_0 \theta$ sufficiently small 
$$\mu_k^h\to\mu^h\quad\text{in the varifold sense},$$
where $\mu^h$ is an integral, rectifiable 2-varifold with weak mean curvature vector $H_h\in L^2(\mu^h)$ such that by lower semicontinuity
$$W_h(\mu^h)=\frac{1}{4}\int|H_h|^2 d\mu^h \leq \liminf_{k\to\infty} W_h(\Sigma_k)=\inf_{[\Sp^2,(\Rtre,\delta+h)]}W_h<4\pi.$$
Now our candidate for a minimizer is given by
$$\Sigma=\spt\mu^h.$$
Now it follows from the monotonicity formula as in \cite{SiL} that
$$\spt\mu^h_k\to\spt\mu^h=\Sigma\quad\text{in the Hausdorff distance sense}.$$
From this convergence and (\ref{lowdiambound}) it follows that
$$\diam_h (\spt \mu^h)>0.$$
Moreover remember that due to (\ref{eq:EstW1}) we may assume that for some $\delta_0>0$
$$W_e(\Sigma_k)\le8\pi-\delta_0.$$
Now we define the so called bad points with respect to a given $\varepsilon>0$ in the following way: define the Radon measures $\alpha_k$ on $\Rtre$ by
$$\alpha_k=\mu_k^h\llcorner|A_k^h|^2.$$
From (\ref{eq:EstAe2Ah2}) and the Gauss Bonnet Theorem it follows that $\alpha_k(\Rtre)\le C$ is uniformly bounded, therefore there exists a Radon measure $\alpha$ on $\Rtre$ such that (after passing to a subsequence) $\alpha_k\to\alpha$ weakly as Radon measures. It follows that $\spt\alpha\subset\Sigma$ and $\alpha(\Rtre)\le C$. Now we define the bad points with respect to $\varepsilon>0$ by
\begin{equation}\label{7}
\B=\left\{\xi\in\Sigma\,\big|\,\alpha(\{\xi\})>\varepsilon^2\right\}.
\end{equation}
Since $\alpha(\Rtre)\le C$, there exist only finitely many bad points. Moreover for $\xi_0\in\Sigma\setminus\B$ there exists a $\rho_0=\rho_0(\xi_0,\varepsilon)>0$ such that $\alpha(B^e_{\rho_0}(\xi_0))<\frac{3}{2}\varepsilon^2$, and since $\alpha_k\to\alpha$ weakly as measures we get
\begin{equation}\label{9'}
\int_{B^e_{\rho_0}(\xi_0)}|A^h_k|^2\,d\mu_k^h\le \frac{3}{2}\varepsilon^2\quad\text{for }k\text{ sufficiently large.}
\end{equation}
%where, as before, $A^h_k$ and ${\cal H}^2_h$ denote the second fundamental form of $\Sigma_k$ and the 2-dimensional Hausdorff measure in $(\Rtre,\delta+h)$.
Consider geodesic normal  coordinates of the Riemannian manifold $(\Rtre,\delta+h)$ centered at $\xi_0$ (the coordinates of $\xi_0$ are $0$); in these coordinates the metric can be written as (see for example \cite{LP} formula (5.4) page 61)
\begin{equation}\label{eq:NormCoord}
(\delta+h)_{\mu\nu}(x)= \d_{\mu\nu} + \frac{1}{3} R_{\mu \sigma \lambda \nu} x^\sigma x^\lambda+O(|x|^3)
             =            \d_{\mu\nu} + o_1(1)(x)_{\mu\nu}, 
\end{equation}
where as before $|o_1(1)(x)|+|\D o_1(1)(x)|\to 0$  for $x\to 0$. Called $inj(\xi_0)>0$ the injectivity radius at $\xi_0$, for $\rh_0<inj(\xi_0)$ we can put on $B_{\rho_0}(\xi_0)$ the normal coordinates just introduced and work on $\Sigma_k\cap B_{\rho_0}(\xi_0)$ as it was immersed in the manifold $(\Rtre, \delta+\tilde{h})$, where $\|\tilde h\|_{C^1}$ can be taken arbitrarily small (for $\rh_0$ small enough). Then taking $\gamma>0 $ sufficiently small in estimate \eqref{eq:EstAe2Ah2}, using \eqref{eq:EstAreaForm} and Proposition \ref{prop:MFh}, we conclude that for $\rh_0$ small enough the bound \eqref{9'} implies

\begin{equation}\label{9}
\int_{\Sigma_k\cap B^e_{\rho_0}(\xi_0)}|A^e_k|^2\,d{\cal H}^2_e\le 2\varepsilon^2\quad\text{for }k\text{ sufficiently large}.
\end{equation}
Now fix $\xi_0\in\Sigma\setminus\B$ and let $\rho_0$ as in (\ref{9}). Let $\xi\in\Sigma\cap B_\frac{\rho_0}{2}(\xi_0)$. We want to apply Simon's graphical decomposition lemma to show that the surfaces $\Sigma_k$ can be written as a graph with small Lipschitz norm together with some "pimples" with small diameter in a neighborhood around the point $\xi$. This is done in exactly the same way Simon did in \cite{SiL}. We just sketch this procedure. By the Hausdorff convergence there exists a sequence $\xi_k\in\Sigma_k$ such that $\xi_k\to\xi$. In view of (\ref{9}) and the monotonicity formula applied to $\Sigma_k$ and $\xi_k$ the assumptions of Simon's graphical decomposition lemma are satisfied for $\rho\le\frac{\rho_0}{4}$ and infinitely many $k\in\N$. Since $W_e(\Sigma_k)\le8\pi-\delta_0$, we can apply Lemma 1.4 in \cite{SiL} to deduce that for $\theta\in\left(0,\frac{1}{2}\right)$ small enough, $\tau\in\left(\frac{\rho}{4},\frac{\rho}{2}\right)$ and infinitely many $k\in\N$ only one of the discs $D_{\tau,l}^k$ appearing in the graphical decomposition lemma can intersect the ball $B_{\theta\frac{\rho}{4}}(\xi_k)$ at fixed $k$. Moreover, by a slight perturbation from $\xi_k$ to $\xi$, we may assume that $\xi\in L_k$ for all $k\in\N$. Now $L_k\to L$ in $\xi+G_2(\Rtre)$, and therefore we may furthermore assume that the planes, on which the graph functions are defined, do not depend on $k\in\N$. After all we get a graphical decomposition in the following way.

\begin{lem}\label{final}
Let $\xi_0\in\Sigma\setminus\B$ and $\rho_0$ as in (\ref{9}). Let $\xi\in\Sigma\cap B_\frac{\rho_0}{2}(\xi_0)$. Then for $\varepsilon\le\varepsilon_0$, $\rho\le\frac{\rho_0}{4}$ and infinitely many $k\in\N$ there exist pairwise disjoint closed subsets $P_1^k,\ldots,P_{N_k}^k$ of $\Sigma_k$ such that 
$$\Sigma_k\cap\overline{B_{\theta\frac{\rho}{8}}(\xi)}=D_k\cap\overline{B_{\theta\frac{\rho}{8}}(\xi)}=\left(\graph u_k\cup\bigcup_n P_n^k\right)\cap\overline{B_{\theta\frac{\rho}{8}}(\xi)},$$
where $D_k$ is a topological disc and where the following holds:  
\begin{enumerate}
\item The sets $P_n^k$ are topological discs disjoint from $\graph u_k$.
\item $u_k\in C^\infty(\overline{\Omega_k},L^\perp)$, where $L\subset\Rtre$ is a 2-dim. plane with $\xi\in L$, and $\Omega_k=\left(B_{\lambda_k}(\xi)\cap L\right)\setminus\bigcup_m d_{k,m}$. Here $\lambda_k>\frac{\rho}{4}$ and the sets $d_{k,m}\subset L$ are pairwise disjoint closed discs.
\item The following inequalities hold:
\begin{eqnarray}
& & \hspace{-2cm}\sum_m\diam d_{k,m}+\sum_n\diam P_n^k\le c\left(\int_{\Sigma_k\cap B^e_{2\rho}(\xi)}|A^e_k|^2\,d{\cal H}^2_e\right)^\frac{1}{4}\rho\le c\varepsilon^\frac{1}{2}\rho, \\
& & ||u_k||_{L^\infty(\Omega_k)}\le c\varepsilon^{\frac{1}{6}}\rho+\delta_k,\quad\text{where }\delta_k\to0,\phantom{\left(\int_{B_{2\rho}}\right)^\frac{1}{4}} \\
& & ||\D u_k||_{L^\infty(\Omega_k)}\le c\varepsilon^{\frac{1}{6}}+\delta_k,\quad\text{where }\delta_k\to0.\phantom{\left(\int_{B_{2\rho}}\right)^\frac{1}{4}}
\end{eqnarray}
\end{enumerate} 
\end{lem}
%\subsection{$C^\infty$ regularity of  $\Sigma$ }\label{Subsec:CinftyRegulatity}
%Since this semiperturbative setting is closely related with the setting in \cite{SiL}, we just sketch the procedure for proving regularity pointing out the main differences with \cite{SiL} and refering to the mentioned paper for more  details (for details see also \cite{MSComp} or \cite{Schy}).
In the next step one proves a power decay for the $L^2$-norm of the second fundamental form on small balls around the good points $\xi\in\Sigma\setminus\B$. This will help us to show that $\Sigma$ is actually $C^{1,\alpha}\cap W^{2,2}$ away from the bad points.

\begin{lem}\label{2ff-absch}
Let $\xi_0\in\Sigma\setminus\B$. There exists a $\rho_0=\rho_0(\xi_0,\varepsilon)>0$ such that for all $\xi\in\Sigma\cap B_{\frac{\rho_0}{2}}(\xi_0)$ and all $\rho\le\frac{\rho_0}{4}$ we have
$$\liminf_{k\to\infty}\int_{\Sigma_k\cap B^e_{\theta\frac{\rho}{8}}(\xi)}|A^e_k|^2\,d{\cal H}^2_e \le c\rho^\alpha,$$
where $\alpha\in(0,1)$ and $c<\infty$ are universal constants.
\end{lem}
The proof of this Lemma is the same as in \cite{SiL}, noticing that in view of the expansion of the metric in normal coordinates as above one can pass from the setting $(\Rtre,\delta+h)$ to the standard euclidean setting up to an error bounded by $c\rho^2$ (for more details see also the proof Lemma 3.6 in \cite{KMS}).\\
\\
Next one shows that the candidate minimizer $\Sigma$ is given locally by a Lipschitz graph with small Lipschitz norm away from the bad points. Again we briefly sketch the construction, for more details see the aforementioned papers. First of all one replaces the pimples of the Graphical Decomposition Lemma \ref{final} with appropriate  graph extensions with small $C^1$ norm, thus they converge to a Lipschitz function with small Lipschitz norm. Then, using a generalized Poincar\'e inequality proved in Lemma A.1 in \cite{SiL} together with the previous Lemma \ref{2ff-absch}, one proves that for all $\xi\in \Sigma \cap B^e_{\frac{\rho_0}{2}}(\xi_0)$ and all sufficiently small $\rh$ 
\begin{equation}\label{mu=graph}
\mu^h\llcorner B^e_\rho(\xi)={\cal H}^2_h \llcorner \left(\graph u\cap B^e_\rho(\xi)\right),
\end{equation} 
where $u\in C^{0,1}(B^e_\rho(\xi)\cap L,L^\perp)$. For more details see the proof of Lemma 3.7 in \cite{KMS}.\\
\\
Since the limit measure $\mu^h$ has weak mean curvature $H_h\in L^2(\mu^h)$, it follows from the definition of the weak mean curvature that $u \in W^{2,2}$; moreover using Lemma \ref{2ff-absch} one can show that the $L^2$ norm of the Hessian of $u$ satisfies the following power decay
\begin{equation}\label{30}
\int_{B_\rho\cap L}|\D^2 u|^2\le c\rho^\alpha.
\end{equation}
From Morrey's lemma (see \cite{GT}, Theorem 7.19) it follows that $u\in C^{1,\alpha}\cap W^{2,2}$. Thus our candidate minimizer can be written as a $C^{1,\alpha}\cap W^{2,2}$-graph away from the bad points. \\
\\
Now one excludes the bad points $\B$ by proving a similar power decay as in Lemma \ref{2ff-absch} for balls around the bad points. This relies on the fact that we are minimizing among spheres. For details see \cite{KMS}, pages 17ff. Therefore our candidate minimizer is given locally by a $C^{1,\alpha}\cap W^{2,2}$-graph everywhere.\\
\\
Again as in \cite{KMS} one can now show that $\Sigma$ is actually a topological sphere.
%Let us point out that by \cite{SiL}, $\gen (\Sigma)\leq \liminf_k \gen(\Sigma_k)=0$ (for a different proof see Lemma 6.2 in \cite{MSComp}). 
Via a standard approximation argument one can check that  
$$\inf \{W_h(\Sigma)| \Sigma \text{ is a smooth embedded $2$-sphere}\}= \inf\{W_h(\Sigma)|  \Sigma \text{ is a } C^1\cap W^{2,2}\text{-embedded $2$-sphere}\}.$$ 
Then by lower semicontinuity of the Willmore energy as mentioned before and the strict $8\pi$ bound of the euclidean Willmore energy it follows that $\Sigma$ is an embedded $2$-sphere which minimizes $W_h$ among $C^{1}\cap W^{2,2}$-embedded $2$-spheres, in particular it satisfies a fourth order Euler Lagrange equation, which fits into the scheme of Lemma 3.2 in \cite{SiL}. Higher regularity and actually smoothness follows as in \cite{SiL}, for more details see again \cite{KMS}. Therefore Theorem \ref{thm:ExW} is proved.
%\begin{equation*}\label{eq:Wh'}
%W_h'(\Sigma)=\frac{1}{2} \triangle H-\frac{1}{4}H(H^2-2|A|^2-2Ric_h(\nu,\nu))
%\end{equation*}
%where $\triangle$ is the Laplace Beltrami of the surface $\Sigma$ and $Ric_h(\nu,\nu)$ is the Ricci tensor of $(\Rtre, \delta+h)$ evaluated on the unit normal $\nu$ to $\Sigma$.
%It is a long and tedious computation but it is possible to check that the Euler Lagrange equation of $W_h$ fits in Lemma 3.2 in \cite{SiL}. 

%It follows that the function $u$ locally representing $\mu^h$  is actually $C^{2,\alpha}\cap W^{3,2}$ and the $L^2$ norm of the $3^{rd}$ derivatives satisfies the  power decay $\int_{B_\rh} |\D^3 u|^2 \leq c\rho^{\alpha}$.
%Now  using the  difference quotients method one proves that the function $u$ is actually $C^{3,\alpha}\cap W^{4,2}$ and the $L^2$ norm of the $4^{th}$ derivatives satisfies the  power decay  $\int_{B_\rh} |\D^4 u|^2 \leq c\rho^{\alpha}$; continuing this bootstrap argument one shows the  smoothness of $u$ and thus of $\Sigma$.   

\section{Proof of Theorem \ref{thm:ExE1} and Theorem \ref{thm:ExW1}}

In this section we prove Theorems \ref{thm:ExE1} and  \ref{thm:ExW1}. Recall the assumptions on the ambient manifold:  $(M,h)$ is a non compact 3-manifold without boundary, of bounded geometry (i.e. satisfying \eqref{Cond:BSC} and \eqref{Cond:InjRad} ) which is either asymptotically euclidean as in $iiia)$ of the Introduction or is hyperbolic outside a compact subset as in $iiib)$.  
%Let $(M,h)$ be a non-compact 3-dimensional Riemannian manifold without boundary, with bounded sectional curvature $|K^M| \leq \Lambda<\infty$ and strictly positive injectivity radius $Inj(M,h)\geq \bar{\rh} >0$. We assume that 
%\\- either the manifold is hyperbolic outside a compact subset, namely there exists $\Omega \subset \subset M$ such that the sectional curvature $K^M\leq 0$ on $M\setminus \Omega$,  
%\\- or $(M,h)$ is asymptotically euclidean in the following sense: there exist compact subsets $K_1 \subset \subset M $ and $K_2 \subset \subset \Rtre$ such that
%\begin{equation}
%(M \backslash K_1) \text{ is isometric to } (\Rtre\backslash K_2, \d+o_1(1)),
%\end{equation}
%where $(\Rtre, \d+o_1(1))$ denotes the Riemannian manifold $\Rtre$ endowed with the euclidean metric $\d_{ij}+o_1(1)_{ij}$ and $o_1(1)$ denotes a symmetric bilinear form which goes to 0 with its first derivatives at infinity, namely
%$$\lim_{|x|\to \infty} (|o_1(1)(x)|+|\nabla o_1(1)(x)|)=0.$$

\subsection{A priori estimates for a minimizing sequence of $E_1$ and $W_1$}
In this subsection we prove the geometric a priori estimates on minimizing sequences of $E_1$ and $W_1$ needed for having compactness and non degeneracy; namely we prove lower and upper bounds on the diameters and we show that the minimizing sequences cannot escape to infinity (the upper bound on the area  clearly follows from the expression of $W_1,E_1$). 
Since the ambient manifold is non compact, it is not trivial a priori that the minimizing sequences have a uniform upper diameter bound. But actually this holds, and it is proved below after a local monotonicity formula (a similar monotonicity formula has been obtained independently by Link in his Ph.D. Thesis, see \cite{FL}).

\begin{lem}\label{lem:LocMF}
Let $(M,h)$ be a (maybe non compact) 3-manifold of bounded geometry, i.e.  satisfying \eqref{Cond:BSC} and \eqref{Cond:InjRad}.
Consider a smooth surface $\Sigma$ immersed in $(M,h)$ and fix $x_0 \in M$. Then there exists a radius $\rho_0=\rho_0(\Lambda,\bar{\rh})$ and constant $C_{\Lambda,\bar{\rh}}$ depending just on the bounds on the injectivity radius and the sectional curvature but independent of $x_0$ such that for any $0<\sigma<\rho<\rho_0$ the following local monotonicity formula holds:

\begin{equation}\label{eq:LocalMF}
\sigma^{-2}|\Sigma \cap B_{\sigma}(x_0)|\leq C_{\Lambda,\bar{\rh}} \left(\rh^{-2} |\Sigma \cap B_{\rh}(x_0)|+ E(\Sigma \cap B_{\rh}(x_0)) \right),
\end{equation}

where $E(\Sigma \cap B_{\rh}(x_0)):= \frac{1}{2} \int_{\Sigma \cap B_{\rh}(x_0)} |A|^2 d\mu_g$.
\end{lem}

\begin{pf}
Fix a point $x_0 \in M$ and on the metric ball $B_{\bar{\rh}} (x_0)\subset M$ consider Riemann normal coordinates centered in $x_0$, i.e. $x_0$ is the origin in the coordinate system. As explained before in \eqref{eq:NormCoord}, in these coordinates the metric $h$ of $M$ is a perturbation of the euclidean metric in the coordinate system:
$$h_{\mu \nu}(x)=\delta_{\mu \nu}+ \tilde{h}^{x_0}_{\mu\nu}(x)=\delta_{\mu\nu}+o^{x_0}_1(1)(x)_{\mu\nu},$$
where the remainder $|o^{x_0}_1(1)(x)|+|\D o^{x_0}_1(1)(x)|\to 0$  for $x\to 0$ uniformly with respect to 
$x_0$ thanks to assumptions \eqref{Cond:BSC} and \eqref{Cond:InjRad}. Let us recall the euclidean monotonicity formula of Simon (formula 1.3 in \cite{SiL}):
\begin{equation}\label{EuMF}
(2\sigma)^{-2} |\Sigma \cap B^e_{2\sigma}(x_0)|_e\leq C\left((\rh/2)^{-2}|\Sigma \cap B^e_{\rh/2}(x_0)|_e+W_e(\Sigma \cap B^e_{\rh/2}(x_0))\right)
\end{equation}
for $0<2\sigma<\rh/2<\bar{\rh}$. For $0<\sigma<2\sigma<\rh/2<\rh<\rho_0=\rho_0(\bar {\rh}, \Lambda)$ small enough, using Lemma \ref{lem:Dist} and Lemma \ref{lem:AreaEst} we estimate the area term as follows
$$\frac{1}{C} |\Sigma \cap B^h_{\sigma}(x_0)|_h\leq  |\Sigma \cap B^e_{2\sigma}(x_0)|_e ,\quad  |\Sigma \cap B^e_{\rh/2}(x_0)|_e\leq C |\Sigma \cap B^h_{\rh}(x_0)|_h, $$
where $B^e_{\sigma}(x_0)$ and $B^h_{\sigma}(x_0)$ are the balls in the coordinate metric and in metric $h$, $|\,.\,|_e$ and $|\,.\,|_h$  are the areas in the coordinate metric and in metric $h$. Now let us bound the Willmore term. Using Lemma \ref{lem:estH} and estimate \eqref{eq:EstAe2Ah2}, since $|H_h|^2\leq 2 |A_h|^2$ we get  
$$|H_e|^2\leq C (|H_h|^2+|A_h|^2+1) \leq C(|A_h|^2+1)\quad \text{in } B_{\rh_0}(x_0) \text{ for $\rho_0=\rho_0(\bar {\rh}, \Lambda)$ small enough,}$$
which integrated gives (we use again Lemma \ref{lem:Dist} and Lemma \ref{lem:AreaEst})
\begin{equation}\label{eq:WeWh}
W_e(\Sigma \cap B^e_{\rh/2}(x_0))\leq C\left( \int_{\Sigma \cap B^h_{\rh}(x_0)}|A_h|^2 \sqrt{\mathring{\delta+h}}+|\Sigma \cap B^h_{\rh}(x_0)|_h\right).
\end{equation}
We conclude that
$$\sigma^{-2}|\Sigma \cap B^h_{\sigma}(x_0)|_h\leq C_{\Lambda,\bar{\rh}} \left(\rh^{-2} |\Sigma \cap B^h_{\rh}(x_0)|_h +\int_{\Sigma \cap B^h_{\rh}(x_0)}|A_h|^2 \sqrt{\mathring{\delta+h}} \right)$$

for a constant $C_{\Lambda,\bar{\rh}}$ depending just on the bounds on the injectivity radius and the sectional curvature but independent on the base point $x_0$.
\end{pf}

\begin{pro}\label{pro:UpDiamWe}
Let $(M,h)$ be a (maybe non compact) Riemannian 3-manifold of bounded geometry, i.e. satisfying \eqref{Cond:BSC} and \eqref{Cond:InjRad}. 

Then there exists a constant $C=C(\bar{\rh}, \Lambda)>0$ such that for every connected, smooth, closed, immersed, oriented surface $\Sigma \hookrightarrow (M,h)$ we have
$$\diam \Sig \leq  \max\{1,C (\mu_g(\Sigma)+W(\Sigma)-\chi_E(\Sigma))\},$$
where $\mu_g(\Sigma)$, $W(\Sigma)$ and $\chi_E(\Sigma)$ are the area, the Willmore functional and the Euler characteristic of $\Sigma$. 
\end{pro}

\begin{pf}
We may assume that $\diam_g \Sig \geq 1$, otherwise the proposition follows immediately. Since $(M,h)$ is of bounded geometry, by Lemma \ref{lem:LocMF} there exists a constant $C=C(\bar{\rh},\Lambda)$ such that for $0<\sigma < \rh < \rh_0 = \rho_0(\bar{\rho},\Lambda)$ the local monotonicity formula \eqref{eq:LocalMF} holds, namely
$$\sigma^{-2}\mu_g(\Sigma\cap B_\sigma(x)) \leq C(\rho^{-2}\mu_g(\Sigma\cap B_\rh(x)) + E(\Sigma \cap B_\rh(x))).$$
Letting $\sigma\to 0$ it follows for every $\rh\leq \rh_0$ and $x\in \Sigma$ that
\begin{equation}\label{eq:AreaBel}
1\leq C(\rho^{-2}\mu_g(\Sigma\cap B_\rh(x)) + E(\Sigma \cap B_\rh(x))).
\end{equation}
Since $\Sigma$ is compact, there exists a pair of points $x,y \in \Sigma$ such that $d(x,y)=\diam\Sigma$. Let 
$$\frac{1}{2} \min (1, \rh_0) < \rh < \min (1, \rh_0)< \diam\Sigma.$$
Let $N\geq 1$ be such that $\frac{1}{2}\diam\Sigma\le N\rho\le\diam\Sigma$ and define for $i=1,\ldots,N$ the sets
%Consider the corresponding partition of the metric ball $B_{\diam_g\Sig}(x)$ into $N$ spherical (metric) annuli at distance $\rh$ one to each other
$$A_i=B_{i\rh}(x)\backslash B_{(i-1)\rh}(x),$$
where $B_{i\rh}(x)$ is the metric ball. Since the surface $\Sigma$ is connected, for each annulus $A_i$ there exists a metric ball $B_{\frac{\rh}{3}}(x_i)\subset A_i$ with center $x_i \in \Sigma$. For each ball $B_{\frac{\rh}{3}}(x_i)$ we can apply the estimate \eqref{eq:AreaBel}. Since the balls $B_{\frac{\rh}{3}}(x_i)$ are pairwise disjoint, summing over $i$ yields
$$N\le C \;\sum_{i=1}^N(\rho^{-2}\mu_g(\Sigma\cap B_{\frac{\rh}{3}}(x_i)) + E(\Sigma\cap B_{\frac{\rh}{3}}(x_i) ))\le C \;(\rho^{-2}\mu_g(\Sigma) + E(\Sigma)).$$
%where the last inequality comes from the disjointness of the balls $B_{\frac{\rh}{3}}(x_i)$. 
Multiplying both sides by $\rh^2$ it follows since $\rh\leq 1$ that
$$\rh\diam\Sigma\le 2N\rh^2\leq C (\mu_g(\Sigma)+ E(\Sigma)).$$
%where the last inequality comes from the condition $\rh\leq 1$. 
By definition of $\rho$ we have $\frac{1}{\rh}< 2 \max (1,1/\rh_0)\leq C=C(\bar{\rh},\Lambda)$, so
\begin{equation}\label{eq:diamProv}
\diam\Sigma\leq C (\mu_g(\Sigma)+ E(\Sigma)).
\end{equation}
Now, by the Gauss equation,  observe that
\begin{equation} \label{eq:G}
\frac{1}{4}|H|^2 - \frac{1}{2}|A^\circ|^2
= \frac{1}{2} \big(|H|^2 - |A|^2\big) = K_g-K^M,
\end{equation}
where $K_g$ is the sectional curvature (also called Gauss curvature) of the induced metric on $\Sigma$ and $K^M$ is the sectional curvature of the tangent plane of $\Sigma$ in $TM$. Integrating \eqref{eq:G}, by Gauss Bonnet Theorem we obtain
\begin{equation}\label{eq:E<W}
E(\Sigma):=\frac{1}{2}\int_{\Sigma}|A|^2 d\mu_g\leq \frac{1}{2}\int_{\Sigma}|H|^2 d\mu_g+\Lambda\;\mu_g(\Sigma) -2\pi\chi_E(\Sigma)= 2 W(\Sigma)+\Lambda \; \mu_g(\Sigma)-2\pi\chi_E(\Sigma)
\end{equation}
and therefore the proposition follows combining \eqref{eq:E<W} and \eqref{eq:diamProv}. 
%where $C_{\bar{\rh},\Lambda}$ is a constant depending on $\bar{\rh}\leq Inj(M)$ and $\Lambda$. We can conclude that
%$$ \diam_g \Sigma \leq C_{\bar{\rh},\Lambda} \; (|\Sigma|_g+ W(\Sigma)).$$ 
\end{pf}

In order to prove an upper and lower bound on the  diameters, we first show that the infimum of $W_1$ and $E_1$ is strictly less than $4\pi$, assuming that there exists a point $\bar{x}\in M$ where the scalar curvature is greater than $6$.
%In order to get this bound it is useful to prove that the $\inf W_1$ in the manifold is less then the corresponding $\inf$ in euclidean space; so let us compute the expansion of $W_1$ on small geodesic spheres.

\begin{lem} \label{lem:ExpW1}
Let $(M,h)$ be a (maybe non-compact) Riemannian 3-manifold. Assume there exists a point $\bar{x}\in M$ where the scalar curvature is greater than $6$, namely 
$$R^M(\bar{x})>6.$$
Then there exist $\epsilon>0$ and $\rh>0$ such that the geodesic sphere $S_{\bar{x},\rh}$ of center $\bar{x}$ and radius $\rh$ satisfies
$$E_1(S_{\bar{x},\rh})=\int_{S_{\bar{x},\rh}}\left(\frac{|A|^2}{2}+1\right)  d\mu_g<4\pi-2\epsilon, $$
$$W_1(S_{\bar{x},\rh})= \int_{S_{\bar{x},\rh}}\left(\frac{|H|^2}{4}+1\right) d\mu_g<4\pi-2\epsilon.$$
\end{lem}

\begin{pf}
From Proposition 3.1 of \cite{Mon1} it follows that on the geodesic spheres $S_{\bar{x},\rh}$ one has
$$W(S_{\bar{x},\rh})=\frac{1}{4} \int_{S_{\bar{x},\rh}} |H|^2 d\mu_g=4\pi-\frac{2\pi}{3}R^M(\bar{x})\rh^2+O(\rh^3).$$
From equation (8) in the proof of Proposition 3.1 in \cite{Mon1} it follows that
$$|S_{\bar{x},\rh}|_g=4 \pi \rh^2+O(\rh^4).$$
Hence the expansion of $W_1$ on small geodesic spheres is
$$W_1(S_{\bar{x},\rh})= 4\pi-\Big(\frac{2}{3}R^M(\bar{x})- 4 \Big) \pi  \rh^2+O(\rh^3).$$
Thus if $R^M (\bar{x})>6$, for $\rh>0$ and $\e>0$ sufficiently small the second inequality follows.\\
\\
For the first inequality observe that $\frac{1}{2}|A|^2=\frac{1}{4}|H|^2+\frac{1}{2}|A^\circ|^2$. Moreover
$$\frac{1}{2} \int_{S_{\bar{x},\rh}} |A ^\circ|^2 d\mu_g=\frac{1}{4}\int_{S_{\bar{x},\rh}} (k_1-k_2)^2 d\mu_g=\int_{S_{\bar{x},\rh}} \Big(\frac{|H|^2}{4}-k_1k_2\Big)d\mu_g$$
is the so called Conformal Willmore functional and was studied by the first author in \cite{Mon2}. In the cited paper the expansion of the functional on geodesic spheres of small radius is computed, and it follows by putting $w=0$ 
%($w$ is the perturbation of the geodesic sphere) 
in Lemma 3.5 and Proposition 3.8 of \cite{Mon2} that
$$\frac{1}{2} \int_{S_{\bar{x},\rh}} |A ^\circ|^2 d\mu_g= O(\rh ^4).$$
Therefore $E_1(S_{\bar{x},\rh})=W_1(S_{\bar{x},\rh})+O(\rh^4)$, and the first inequality follows from the second one. 
\end{pf}

Thanks to Proposition \ref{pro:UpDiamWe}, Lemma \ref{lem:ExpW1} and Remark \ref{rem:WeEucl}, we show in the next step that minimizing sequences for the functional $E_1$, respectively $W_1$, stay in a compact subset of the manifold $M$.
    
%\begin{pro}\label{pro:fkCompact}
%Let $(M,g)$ be a non compact  Riemannian $3$-manifold  with bounded geometry (i.e. with bounded sectional curvature $|\bar{K}| \leq \Lambda^2$ for  some $\Lambda \in \R,$ and striclty positive injectivity radius $Inj(M,g)\geq \bar{\rh} >0$) such that
%\\i) or $(M,g)$ is asymptotically euclidean in the sense of  \eqref{eq:defAsFlat} or outside a compact subset has negative sectional curvature: $\exists \Omega \subset \subset M$ such that $\bar K\leq 0$ on $M\setminus \Omega$
%\\ii) the  scalar curvature is striclty greater than $6$ at a point:
%$$\exists \bar{p} \in M: \quad R_g(\bar{p})>6.$$
%Consider a minimizing sequence $\{f_k:\Sp^2 \hookrightarrow M\}_{k\in \N}$ of smooth immersions of $W_1$ (respectively of $E_1$) among immersions of the same kind.     
%\\Then there exists a compact subset $K\subset \subset M$ such that $f_k(\Sp^2)\subset K$ for all $k\in \N$.
%\end{pro}    

\begin{pro}\label{pro:fkCompact}
Let $(M,h)$ be a non compact Riemannian $3$-manifold without boundary with bounded geometry, i.e. satisfying \eqref{Cond:BSC} and \eqref{Cond:InjRad}, with asymptotic behavior as in $iiia)$ or in $iiib)$. Assume that the scalar curvature is strictly greater than $6$ at a point $\bar{x} \in M$, namely
$$R^M(\bar{x})>6.$$

Let $f_k:\Sp^2 \hookrightarrow M$ be a minimizing sequence for $E_1$, respectively $W_1$. Then there exists a compact subset $K\subset \subset M$ such that $f_k(\Sp^2)\subset K$ for all $k\in \N$.
\end{pro}    

\begin{pf}
From the assumption on the scalar curvature it follows from Lemma \ref{lem:ExpW1} that
\begin{equation}\label{eq:W1<4pi}
\lim_{k\to\infty} E_1(f_k)\leq 4\pi -2\e,\text{ respectively }\lim_{k\to\infty} W_1(f_k)\leq 4\pi-2\e.
\end{equation}
Since $\frac{1}{2}|A|^2=\frac{1}{4}|H|^2+\frac{1}{2}|A^\circ|^2$, clearly 
\be \label{eq:WW1E1}
 W(f)\leq W_1(f)\le E_1(f) \quad \forall f \in [\Sp^2,M], 
\ee 
and Proposition \ref{pro:UpDiamWe} implies
\be \label{eq:udb}
\diam_M (f_k(\Sp^2))\leq 1+C (\mu_{g_k}(\Sp^2)+W(f_k))\leq C
\ee
for some constant $C<\infty$ independent of $k$.
%Observe that if instead $\{f_k\}_{k\in \N}$ is a minimizing sequence of $E_1$, Lemma \ref{lem:ExpW1} implies $\lim_k E_1(f_k)\leq 4\pi -2\e $, but as shown in the proof of the previously cited Lemma, $E_1(f_k)\geq W_1(f_k)$, so inequality \eqref{eq:W1<4pi} also holds for minimizing sequences of $E_1$.

Let us first consider the case $(M,h)$ asymptotically euclidean and $f_k$ minimizing sequence for $W_1$; if the thesis is not true, then, up to subsequences, for every $k \in \N$ we can take a point $\xi_k \in f_k(\Sp^2)\subset \Rtre$ (recall that outside a compact subset, $(M,h)$ is isometric to $(\Rtre, \d+o_1(1)$) such that $|\xi_k| \to \infty$. Since by \eqref{eq:udb} we have that $\diam f_k(\Sp^2)\leq C$ %and since by the minimizing sequence property we trivially have $\mu_{g_k}(\Sp^2)\leq C$ 
,it follows that
$$\liminf_{k\to\infty} \|o_1(1)\|_{C^1(f_k(\Sp^2))}=0. $$
%Now for every fixed $k \in \N$ we can put $\gamma=r_0=\diam_e (f_k(\Sp^2))\stackrel{k\to \infty} \rightarrow 0$ in estimate \eqref{eq:EstW} (repeat the proof of  Lemma \ref{lem:EstW} with such quantities), since for $k$ large also $\eta$ and $\theta$ are arbitrarily small, passing to the liminf in estimate \eqref{eq:EstW} we can conclude that
Repeating the proof of Lemma \ref{lem:EstW} yields
$$\liminf_{k\to\infty} W_1(f_k)\ge \liminf_{k\to\infty} W(f_k)\geq \liminf_{k\to\infty} W_e(f_k)\geq 4\pi,$$
which contradicts \eqref{eq:W1<4pi}. Thus there exists a compact subset $K\subset \subset M$ such that $f_k(\Sp^2)\subset K$ for all $k \in \N$.  The case $(M,h)$ asymptotically euclidean and $f_k$ minimizing sequence of $E_1$ follows by \eqref{eq:WW1E1}: repeating the arguments above we again arrive to contradict \eqref{eq:W1<4pi}.

Now consider the case $(M,h)$ hyperbolic outside a compact subset: there exists $\Omega \subset \subset M$ such that the sectional curvature $K^M\leq 0$ on $M\setminus \Omega$. The Gauss equation \eqref{eq:G} implies that on $f_k(\Sp^2)\cap (M\setminus \Omega)$ one has
\be\label{eq:HKg}
\frac{1}{4}|H|^2 \geq K_g. 
\ee
If by contradiction the sequence $f_k$ is not contained in any compact subset of $M$, then it follows from the diameter bound that, up to subsequences, $f_k(\Sp^2) \subset M\setminus \Omega$. Since we are working on spheres, integrating \eqref{eq:HKg} and using the Gauss-Bonnet Theorem yields 
\begin{equation}\label{eq:AH>4pi}
\frac{1}{2}\int |A|^2\,d\mu_g\geq \frac{1}{4}\int |H|^2\,d\mu_g \geq 4 \pi,
\end{equation} 
 which implies $E_1(f_k)\geq W_1(f_k)\geq 4\pi$, contradicting \eqref{eq:W1<4pi}. 
\end{pf}

Now we conclude that the minimizing sequences $f_k\in [\Sp^2,M]$ for $W_1$ or $E_1$ cannot shrink to a point, namely
\be\label{prop:LBdiamW1}
\liminf_{k \to \infty} (\diam_M (f_k(\Sp^2)))>0. 
\ee 
Indeed by Proposition \ref{pro:fkCompact} there exists a compact subset $K\subset \subset M$ containing all the surfaces, up to subsequences: $f_k(\Sp^2)\subset K$.  By Lemma \ref{lem:ExpW1} and inequality \eqref{eq:WW1E1} it follows that on the minimizing sequence we have 
\be\label{eq:W<4pi}
\liminf_k W(f_k)\leq 4\pi-2\epsilon.
\ee
Then \eqref{prop:LBdiamW1} follows from \eqref{eq:W<4pi} together with Proposition 2.5 in \cite{KMS}, which also holds for non-compact Riemannian manifolds without boundary if the minimizing sequence stays in a compact subset.

\subsection{Existence and regularity of minimizers for $E_1$, respectively $W_1$}
Let $(M,h)$ be a non compact Riemannian 3-manifold without boundary as in Proposition \ref{pro:fkCompact}. For the problem of minimizing the functional $W_1$, namely for the proof of Theorem \ref{thm:ExW1}, we assume in addition that the sectional curvature $K^M\le2$. Let $f_k\in[\Sp^2,M]$ be a minimizing sequence for the functional $E_1$, respectively $W_1$. It follows from the previous lemmas and propositions that
\begin{itemize}
\item[i)] there exists a constant $C<\infty$ such that $\mu_k(\Sp^2)\le C$, where $\mu_k$ is the induced area measure,
\item[ii)] there exists a compact subset $K\subset \subset M$ such that $f_k(\Sp^2)\subset K$,
\item[iii)] there exists a constant $0<C<\infty$ such that $\frac{1}{C}\leq \diam (f_k(\Sp^2)) \leq C$.
\end{itemize}
Now observe that it follows directly from (\ref{eq:G}) that if $K^M\le2$, then we can estimate the $L^2$-norm of the second fundamental form by the functional $W_1$, namely we have that
$$\frac{1}{2}\int |A|^2\,d\mu_g\le2W_1(f)-4\pi$$
for every immersion $f\in[\Sp^2,M]$. Therefore, no matter if $f_k$ is a minimizing sequence for $E_1$ or $W_1$, it follows in addition from Lemma \ref{lem:ExpW1} that
$$\hspace{-6,2cm}\text{vi)}\hspace{4,5cm}\limsup_{k\to\infty} \frac{1}{2}\int|A_k|^2\,d\mu_k< 4\pi,$$
where $A_k$ denotes the second fundamental form of $f_k$.\\
\\
The properties above are actually all the properties for minimizing sequences for the functional $E_1$, respectively $W_1$, one needs to apply the existence proof in \cite{KMS}. Although in the aforementioned paper it is assumed that the manifold $M$ is compact, we can apply the techniques developed there since minimizing sequences in our setting stay in a compact subset of the non-compact manifold $M$, which is enough. Thus the proof of Theorem 1.1 and Theorem 1.2 in \cite{KMS} can be directly applied in our situation, which proves Theorem \ref{thm:ExE1} and Theorem \ref{thm:ExW1}.

\end{document}